\documentclass[12pt]{amsart}
\usepackage{latexsym}
\usepackage{pstricks}
\usepackage{amsmath}
\usepackage{amssymb}
\usepackage{mathtext}
\newtheorem{propo}{Proposition}[section]

\newtheorem{lemma}[propo]{Lemma}

\newtheorem{theorem}[propo]{Theorem}

\newcommand{\Ker}{\operatorname{Ker}}

\newcommand{\Irr}{{\mathrm {Irr}}}
\newcommand{\IBR}{{\mathrm {IBr}}}

\newcommand{\ord}{{\mathrm {ord}}}

\newcommand{\Ind}{{\mathrm {Ind}}}

\newcommand{\CC}{{\mathbb C}}

\newcommand{\ZZ}{{\mathbb Z}}

\newcommand{\FF}{{\mathbb F}}

\newcommand{\ta}{\hspace{0.5mm}^{2}\hspace*{-0.2mm}}

\begin{document}

\title[irreducible restrictions of  Brauer characters of $G_2(q)$]
{irreducible restrictions of  Brauer characters of the Chevalley
group $G_2(q)$ to its proper subgroups}

\author{Hung Ngoc Nguyen}
\address{Department of Mathematics, University of Florida, Gainesville,
FL 32611, USA} \email{hnguyen@math.ufl.edu}

\begin{abstract}
Let $G_2(q)$ be the Chevalley group of type $G_2$ defined over a
finite field with $q=p^n$ elements, where $p$ is a prime number and
$n$ is a positive integer. In this paper, we determine when the
restriction of an absolutely irreducible representation of $G$ in
characteristic other than $p$ to a maximal subgroup of $G_2(q)$ is
still irreducible. Similar results are obtained for $\ta B_2(q)$ and
$\ta G_2(q)$.
\end{abstract}

\maketitle

\section{Introduction}

Let $G$ be the Chevalley group $G_2(q)$ defined over a finite field
with $q=p^n$ elements, where $p$ is a prime number and $n$ is a
positive integer. Let $M$ be a maximal subgroup of $G$ and $\Phi$ an
absolutely irreducible representation of $G$ in cross characteristic
$\ell$ (i. e. $\ell=0$ or $\ell$ prime, $\ell \nmid q$). Let
$\varphi$ be the irreducible Brauer character of $G$ afforded by
$\Phi$. The purpose of this paper is to find all possibilities of
$\varphi$ and $M$ such that $\varphi|_M$ is also irreducible. This
result is a contribution to the classification of maximal subgroups
of finite classical groups guided by Aschbacher's Theorem \cite{A}.

In order to do that, we need to use the results about maximal
subgroups, character tables, blocks and Brauer trees of $G$ obtained
by many authors.

The list of maximal subgroups of $G$ is determined by Cooperstein in
\cite{Cp} for $p=2$ and Kleidman in \cite{K1} for $p$ odd. The
complex character table of $G$ is determined by Chang and Ree in
\cite{CR} for good primes $p$ (i. e. $p\geq5$). Enomoto computed in
\cite{E} the character table of $G$ when $p=3$. Lastly, when $p=2$,
it is computed by Enomoto and Yamada in \cite{EY}. In a series of
papers \cite{H1}, \cite{HS1}, \cite{HS2}, \cite{S2}, \cite{S3}, and
\cite{S4}, Hiss and Shamash have determined the blocks, Brauer trees
and (almost completely) the decomposition numbers for $G$.

In order to solve the problem, it turns out to be useful to know the
small degrees of irreducible Brauer characters of $G$. The smallest
degrees of non-trivial absolutely irreducible representations of
many quasi-simple groups are collected in \cite{T1} and the result
for $G$ is used in this paper. From the knowledge of the degrees of
irreducible Brauer characters of $G$, we compute the second smallest
such degree. The exact formula for it is given in Theorem \ref{d2l}.
This is also useful for other applications.

Clearly, it suffices to consider maximal subgroups $M$ of $G$ with
$\sqrt{|M|}$ larger than the smallest degree of nontrivial
absolutely irreducible representations of $G$. Thus, we can exclude
many maximal subgroups of $G$ whose orders are small enough by
Reduction Theorem \ref{red}. The remaining maximal subgroups are
treated individually by various methods. The main results of this
paper are the following:

\begin{theorem}\label{main theorem}
Let $G=G_2(q)$, $q=p^n$, $q\geq 5$, $p$ a prime
number. Let $\varphi$ be an absolutely irreducible character of $G$
in cross characteristic $\ell$ and $M$ a maximal subgroup of $G$.
Assume that $\varphi(1)>1$. Then $\varphi|_M$ is irreducible if and
only if one of the following holds:
\begin{enumerate}
  \item[(i)] $q\equiv-1 (\bmod$ $3)$, $M=SL_3(q):2$ and $\varphi$ is the
  unique character of the smallest degree $q^3-1$.
  \item[(ii)] $q\equiv 1 (\bmod$ $3)$, $M=SU_3(q):2$ and $\varphi$ is the
  unique character of the smallest degree. In this case, $\varphi(1)=q^3$
  when $\ell=3$ and $\varphi(1)=q^3+1$ when $\ell\neq3$.
\end{enumerate}
\end{theorem}

The covering groups of $G_2(3)$ and $G_2(4)$ are handled by the
following theorem. Notice that $3\cdot G_2(3)$ has two pairs of
complex conjugate irreducible characters of degree $27$.

\begin{theorem}
\label{main theorem2} Let $G\in \{G_2(3), 3\cdot G_2(3), G_2(4),
2\cdot G_2(4)\}$. Let $\varphi$ be a faithful absolutely irreducible
character of $G$ in cross characteristic $\ell$ of degree larger
than $1$ and $M$ a maximal subgroup of $G$. Then we have:
\begin{enumerate}
  \item[(i)] When $G=G_2(3)$, $\varphi|_M$ is irreducible if and only if one of the following holds:
\begin{enumerate}
  \item[(a)] $M=U_3(3):2$, $\varphi$ is the unique irreducible character of degree
  $14$ when $\ell=7$ or $\varphi$ is any of the irreducible characters of degree
  $14$, $64$ when $\ell\neq 3$, $7$;
  \item[(b)] $M=2^3\cdot L_3(2)$, $\varphi$ is the unique irreducible character
  of degree $14$ when $\ell\neq2$, $3$.
\end{enumerate}
  \item[(ii)] When $G=3\cdot G_2(3)$, the universal cover of $G_2(3)$, $\varphi|_M$
  is irreducible if and only if one of the following holds:
\begin{enumerate}
  \item[(a)] $M=3.P$ or $3.Q$ where $P$, $Q$ are maximal parabolic subgroups of
  $G_2(3)$, $\varphi$ is any of the four irreducible characters of degree $27$;
  \item[(b)] $M=3.(U_3(3):2)$, $\varphi$ is any of the two complex conjugate irreducible
  characters of degree $27$ when $\ell\neq 2$, $3$, $7$;
  \item[(c)] $M=3.(L_3(3):2)$, $\varphi$ is any of the two complex conjugate
  irreducible characters of degree $27$ when $\ell\neq 2$, $3$, $13$;
  \item[(d)] $M=3.(L_2(8):3)$, $\varphi$ is any of the four irreducible characters
  of degree $27$ when $\ell\neq 2$, $3$, $7$.
\end{enumerate}
  \item[(iii)] When $G=G_2(4)$, $\varphi|_M$ is irreducible if and only if one of the following holds:
\begin{enumerate}
  \item[(a)] $M=U_3(4):2$, $\varphi$ is the unique irreducible character of degree $64$
  when $\ell=3$ or $\varphi$ is the unique irreducible character of degree $65$ when $\ell\neq 2$, $3$;
  \item[(b)] $M=J_2$, $\varphi$ is any of the two irreducible characters of
  degree $300$ when $\ell\neq 2$, $3$, $7$.
\end{enumerate}
  \item[(iv)] When $G=2\cdot G_2(4)$, the universal cover of $G_2(4)$,
  $\varphi|_M$ is irreducible if and only if one of the following holds:
\begin{enumerate}
  \item[(a)] $M=2.P$ or $2.Q$ where $P$, $Q$ are maximal parabolic subgroups
  of $G_2(4)$, $\varphi$ is the unique irreducible character of degree
  $12$;
  \item[(b)] $M=2.(U_3(4):2)$, $\varphi$ is the unique irreducible character
  of degree $12$ or $\varphi$ is any of the two irreducible characters of degree
  $104$ when $\ell\neq2, 5$;
  \item[(c)] $M=2.(SL_3(4):2)$, $\varphi$ is the unique irreducible character
  of degree $12$ when $\ell\neq 2, 3$.
\end{enumerate}
\end{enumerate}
\end{theorem}

In the case of complex representations, Theorem \ref{main theorem}
was proved by Saxl [S].

Similar results for the Suzuki group $\ta B_2(q)$ and the Ree group
$\ta G_2(q)$ are obtained in $\S6$.

The paper is divided into six sections. In $\S2$, we set up some
notation and state some lemmata which will be used later. In $\S3$,
we collect the degrees of irreducible Brauer characters and give the
formulas for the first and second degrees of $G$. $\S4$ and $\S5$
are devoted to prove the main results. In the last section, we state
and prove the results for the Suzuki group $\ta B_2(q)$ and the Ree
group $\ta G_2(q)$.

\medskip

{\bf Acknowledgments}. The author is grateful to his advisor,
Professor P. H. Tiep, for devoted guidance and fruitful discussion.
Part of this paper was done while the author was participating in
the Special Semester in Group Representation Theory, Bernoulli
Center, EPFL, Lausanne, Switzerland. The author is grateful to
Professors M. Geck, D. Testerman, and J. Thevenaz as well as the NSA
for financial support.

\section{Preliminaries}

Let $\FF$ be an algebraically closed field of characteristic $\ell$.
Given a finite group $X$, we denote by $\mathfrak{d}_{\ell}(X)$,
$\mathfrak{d}_{2,\ell}(X)$, and $\mathfrak{m}_{\ell}(X)$ the
smallest, the second smallest, and the largest degrees,
respectively, of irreducible $\FF X$-representations of degree
larger than $1$. When $\ell=0$, we use the notation
$\mathfrak{d}_{\mathbb{C}}(X)$, $\mathfrak{d}_{2,\CC}(X)$, and
$\mathfrak{m}_{\CC}(X)$ instead of $\mathfrak{d}_{0}(X)$,
$\mathfrak{d}_{2,0}(X)$, and $\mathfrak{m}_{0}(X)$. If $\chi$ is a
complex character of $X$, we denote by $\widehat{\chi}$ the
restriction of $\chi$ to $\ell$-regular elements of $X$. Throughout
this paper, $\Irr(X)$ is the set of irreducible complex characters
of $X$; $\IBR_\ell(X)$ is the set of irreducible $\ell$-Brauer
characters of $X$; and $\ZZ_n$ is the cyclic group of order $n$.

The following lemmata are well known and we omit most of their
proofs.

\begin{lemma}\label{1}
Let $G$ be a finite group and $H$ be a subgroup of $G$. Let $\Phi$
be an irreducible $\FF G$-representation of degree larger than $1$
such that $\Phi|_H$ is also irreducible. Then
$\mathfrak{m}_{\ell}(H)\geq \deg(\Phi)\geq \mathfrak{d}_{\ell}(G)$.
\end{lemma}

\begin{lemma}\label{2}
Let $G$ be a finite group. Then $\mathfrak{m}_{\ell}(G)\leq
\mathfrak{m}_{\mathbb{C}}(G)$. Furthermore,
$\mathfrak{m}_{\CC}(G)\leq\sqrt{|G/Z(G)|}$.
\end{lemma}

\begin{lemma}\label{lemmaZ}
Let $G$ be a simple group and $V$ an irreducible $\FF G$-module such
that $\dim(V)>1$. Then $Z_G(V):=\{g\in G\mid g|_V=\lambda\cdot
Id_V\text{ for some } \lambda\in\FF\}=1$.
\end{lemma}

\begin{lemma}\label{lemma for Parabolic}
Let $G$ be a finite group and $1\neq A\unlhd H\leq G$. Let $V$ be a
faithful irreducible $\FF G$-module.
\begin{enumerate}
  \item[(i)] Suppose that $C_V(A):=\{v\in V\mid a(v)=v$ for every $a\in A\}\neq 0$.
  Then $V|_H$ is reducible.
  \item[(ii)] Suppose that $O_{\ell}(H)\neq 1$. Then $V|_H$ is reducible.
\end{enumerate}
\end{lemma}

\begin{proof}
(i) Assume the contrary that $V|_H$ is irreducible. By Clifford's
theorem, $V|_A=e\bigoplus_{i=1}^{t}V_i$ where $e$ is the
multiplicity of $V_1$ in $V$ and $\{V_1,...,V_t\}$ is the orbit of
$V_1$ under the action of $H$ on the set of all irreducible $\FF
A$-modules. Since $A$ acts trivially on $C_V(A)\subseteq V$, one of
the $V_i$s is the trivial module. Therefore all $V_1,...,V_t$ are
trivial $A$-modules and $A$ acts trivially on $V$. This contradicts
the faithfulness of $V$.

(ii) Assume that $V|_H$ is irreducible. Then $O_{\ell}(H)$ acts
trivially on $V$. Hence $C_V(O_{\ell}(H))=V$ and we get a
contradiction by part (i).
\end{proof}

\begin{lemma}\label{lemmatau}
Let $\tau$ be an automorphism of a finite group $G$ which fuses two
conjugacy classes of subgroups of $G$ with representatives $H_1$,
$H_2$. Let $\mathcal{A}$ and $\mathcal{B}$ be subsets of
$\IBR_{\ell}(G)$ such that $\tau(\mathcal{A})=\mathcal{B}$,
$\tau(\mathcal{B})=\mathcal{A}$.
\begin{enumerate}
  \item[(i)] Assume $\varphi|_{H_1}$ is irreducible (resp. reducible) for all
$\varphi\in \mathcal{A}\cup\mathcal{B}$. Then $\varphi|_{H_2}$ is
irreducible (resp. reducible) for all $\varphi\in
\mathcal{A}\cup\mathcal{B}$.
  \item[(ii)] Assume $\varphi|_{H_1}$ is irreducible for all
$\varphi\in\mathcal{A}$ and $\varphi|_{H_1}$ is reducible for all
$\varphi\in \mathcal{B}$. Then $\varphi|_{H_2}$ is reducible for all
$\varphi\in\mathcal{A}$ and $\varphi|_{H_2}$ is irreducible for all
$\varphi\in \mathcal{B}$.
\end{enumerate}
\end{lemma}

\begin{proof}
It is enough to show that $\varphi|_{H_1}$ is irreducible (resp.
reducible) for all $\varphi\in \mathcal{A}$ if and only if
$\varphi|_{H_2}$ is irreducible (resp. reducible) for all
$\varphi\in \mathcal{B}$. This is true because for every
$\varphi\in\mathcal{A}$, we have
$\varphi|_{H_1}=\psi\circ\tau|_{\tau(H_2)}=(\psi|_{H_2})\circ\tau$
for some $\psi\in\mathcal{B}$.
\end{proof}

\begin{lemma}
\label{orbitreduction} Let $G$ be a finite group and
$\chi\in\Irr(G)$. Let $H$ be a normal $\ell'$-subgroup of $G$ and
suppose that $\chi|_H=\sum_{i=1}^{t}\theta_i$, where $\theta_1$ is
irreducible and $\theta_1$, $\theta_2$,..., $\theta_t$ are the
distinct $G$-conjugates of $\theta_1$. Then $\widehat{\chi}$ is also
irreducible.
\end{lemma}

\begin{proof}
We have $\widehat{\chi}|_H=\sum_{i=1}^{t}{\theta_i}$ as $H$ is an
$\ell'$-group. Let $\psi$ be an irreducible constituent of
$\widehat{\chi}$. Then there is some $i\in \{1$, $2$, ... , $t\}$
such that ${\theta_i}$ is an irreducible constituent of $\psi|_H$.
Therefore, by Clifford's theorem, all ${\theta_1}$, ${\theta_2}$,
... , ${\theta_t}$ are contained in $\psi|_H$. This implies that
$\widehat{\chi}|_H=\psi|_H$. So $\widehat{\chi}=\psi$ and
$\widehat{\chi}$ is irreducible.
\end{proof}

\begin{lemma}
\label{maximalsubgroup} Let $G$ be a finite group. Suppose that the
universal cover of $G$ is $M.G$ where $M$ is the Schur multiplier of
$G$ of prime order. Then every maximal subgroup of $M.G$ is the
pre-image of a maximal subgroup of $G$ under the natural projection
$\pi: M.G\rightarrow G$.
\end{lemma}

\begin{lemma}\cite{F}
\label{lemmaisaacs} Let $G$ be a finite group and $H\lhd G$. Suppose
$|G:H|=p$ is prime and $\chi\in \IBR_{\ell}(G)$. Then either
\begin{enumerate}
  \item[(i)] $\chi|_H$ is irreducible or
  \item[(ii)] $\chi|_H=\sum_{i=1}^p\theta_i$, where the $\theta_i$s are distinct and irreducible.
\end{enumerate}
\end{lemma}

\begin{lemma}\cite{OT}
\label{Navaro}Let $B$ be an $\ell$-block of group $G$. Assume that
all $\chi\in B\cap \Irr(G)$ are of same degree. Then $B\cap
\IBR_{\ell}(G)=\{\phi\}$ and $\widehat{\chi}=\phi$ for every
$\chi\in B\cap \Irr(G)$.
\end{lemma}

\begin{lemma}\cite[Theorem 1.6]{T2}
\label{lemmaforparabolica1} Let $G$ be a finite group of Lie type,
of simply connected type. Assume that $G$ is not of type $A_1$, $\ta
A_2$, $\ta B_2$, $\ta G_2$, and $B_2$. If $Z$ is a long-root
subgroup and $V$ is a nontrivial irreducible $G$-module, then $Z$
must have nonzero fixed points on $V$.
\end{lemma}

\section{The degrees of irreducible Brauer characters of $G_2(q)$}

In this section, we will recall the value of
$\mathfrak{d}_{\ell}(G)$ and determine the value of
$\mathfrak{d}_{2,\ell}(G)$ when $\ell\nmid q$. The degrees of
irreducible complex characters of $G_2(q)$ can be read off from
\cite{CR}, \cite{E}, \cite{EY} and are listed in Table I.

\medskip
\begin{centering}
{\sc {Table I.}} The degrees of irreducible complex characters of
$G_2(q)$

\begin{center}
\begin{tabular}{|c|c|} \hline
Character  & Degree\\ \hline
 $X_{11}$  & $1$\\\hline
 $X_{12}$&$q^6$\\\hline
$X_{13}, X_{14}$&$\frac{1}{3}q(q^4+q^2+1)$\\\hline
$X_{15}$&$\frac{1}{2}q(q+1)^2(q^2-q+1)$\\ \hline
 $X_{16}$&$\frac{1}{6}q(q+1)^2(q^2+q+1)$\\\hline
 $X_{17}$&$\frac{1}{2}q(q-1)^2(q^2+q+1)$\\\hline
 $X_{18}$&$\frac{1}{6}q(q-1)^2(q^2-q+1)$\\\hline
 $X_{19}, \overline{X_{19}}$&$\frac{1}{3}q(q-1)^2(q+1)^2$\\\hline
 $X_{31}$&$q^3(q^3+\epsilon)$\\\hline
 $X_{32}$&$q^3+\epsilon$\\ \hline
$X_{33}$&$q(q+\epsilon)(q^3+\epsilon)$\\ \hline
$X_{21}$&$q^2(q^4+q^2+1)$\\\hline $X_{22}$&$q^4+q^2+1$\\\hline
$X_{23}, X_{24}$&$q(q^4+q^2+1)$\\\hline $X_{1a}, X_{1b}, X_{2a},
X_{2b}$&$q(q\pm1)(q^4+q^2+1)$\\\hline
 $X'_{1a},X'_{1b}, X'_{2a},X'_{2b}$&$(q\pm1)(q^4+q^2+1)$\\\hline
 $X_{1},X_{2}$ &
$(q\pm1)^2(q^4+q^2+1)$\\\hline
 $X_{a}, X_{b}$& $q^6-1$\\\hline
 $X_{3},X_{6}$ & $(q^2-1)^2(q^2\mp q+1)$\\\hline
\end{tabular}
\end{center}
\end{centering}
\footnotesize where:
\begin{enumerate}
  \item[(i)] $\epsilon\equiv q(\bmod$ $3)$,
  \item[(ii)] $X_{21}$, $X_{22}$, $X_{23}$, $X_{24}$ appear only if $q$ is odd,
  \item[(iii)] $X_{31}$, $X_{32}$, $X_{33}$ appear only if $q$ is not divisible by 3.
\end{enumerate}
\normalsize From the table, we get that if $q\geq5$ then
\begin{equation}\label{d1c}\mathfrak{d}_\mathbb{C}(G)= \left\{\begin {array}{ll}
q^3+1, & q\equiv 1 (\bmod$ $3),\\
q^3-1, & q\equiv 2 (\bmod$ $3),\\
q^4+q^2+1, & q \equiv 0 (\bmod$ $3),
\end {array} \right.
\end{equation} and
\begin{equation}\label{second} \mathfrak{d}_{2,\mathbb{C}}(G)= \left\{\begin
{array}{ll}
\frac{1}{6}q(q-1)^2(q^2-q+1), & p=2$, $3 \text{ or } q=5$, $7,\\
q^4+q^2+1, & p\geq5$, $q>7.
\end {array} \right.
\end{equation}
Moreover, $\Irr(G)$ contains a unique character of degree larger
than $1$ but less than $\mathfrak{d}_{2,\mathbb{C}}(G)$ and this
character has degree $\mathfrak{d}_{\mathbb{C}}(G)$.

The degrees of irreducible $\ell$-Brauer characters of $G_2(q)$ when
$\ell\mid |G|$ and $\ell\nmid q$ can be read off from \cite{H1},
\cite{HS1}, \cite{HS2}, \cite{S2}, and \cite{S3}. They are listed in
Tables II, III for $\ell=2$ and $\ell=3$. When $\ell\geq 5$, we
refer to \cite{H1}, \cite{S2}, and \cite{S3}. In these tables,
$q\equiv \epsilon(\bmod$ $3)$. Moreover, $\varphi_{21}$,
$\varphi_{22}$, $\varphi_{23}$, $\varphi_{24}$ appear only if $q$ is
odd and $\varphi_{31}$, $\varphi_{32}$, $\varphi_{33}$ appear only
if $q$ is not divisible by $3$.

Note that, besides the characters listed in Tables II, III, the
degrees of the remaining irreducible Brauer characters of $G_2(q)$:
$\varphi_{11}$, $\varphi_{17}$, $\varphi_{18}$, $\varphi_{19}$,
$\overline{\varphi_{19}}$, $\varphi'_{1a}$, $\varphi'_{1b}$,
$\varphi'_{2a}$, $\varphi'_{2b}$, $\varphi_{1}$, $\varphi_{2}$,
$\varphi_{a}$, $\varphi_{b}$, $\varphi_{3}$, $\varphi_{6}$ are the
same as the degrees of complex characters with the same indices.

\medskip
\centerline {{\sc {Table II.}} The degrees of irreducible $2$-Brauer
characters of $G_2(q)$ with $q$ odd}
\begin{center}
\begin{tabular}{|c|c|c|} \hline
Character & $4\mid (q-1)$ & $4\mid (q+1)$\\ \hline
 $\varphi_{12}$& \begin{tabular} {c}$\frac{1}{6}(q-1)^2(6q^4+(8-3\alpha-\beta)
 q^3$\\ $+(10-3\alpha+\beta)q^2$\\$+(8-3\alpha-\beta)q+6)$\end{tabular}&
 $\begin{tabular} {c}$\frac{1}{6}(q-1)^2(6q^4+(8-3\alpha-\beta)
 q^3$\\ $+(10-3\alpha+\beta)q^2$\\$+(8-3\alpha-\beta)q+6)$\end{tabular}$\\ \hline
$\varphi_{13},\varphi_{14}$&$\frac{1}{3}(q-1)(q^4+q^3+2q^2+2q+3)$&$\frac{1}{3}(q-1)(q^4+q^3+2q^2+2q+3)$
\\\hline
 $\varphi_{15}$&$q^4+q^2$&$q^4+q^2$\\\hline
$\varphi_{31}$&\begin{tabular} {c}$q^6-1$ if $q\equiv1(\bmod$ $3)$\\
 $(q-1)^2(q^2+1)(q^2+q+1)$ \\if $q\equiv -1(\bmod$ $3)$ \end{tabular}&\begin{tabular} {c}
 $q^6-1$ if $q\equiv1(\bmod$ $3)$\\
 $(q-1)^2(q^2+q+1)$\\$(q^2+(1-\gamma)q+1)$\\ if $q\equiv -1(\bmod$ $3)$
 \end{tabular}\\\hline
 $\varphi_{32}$&$q^3+\epsilon$&$q^3+\epsilon$\\ \hline
$\varphi_{33}$&$q(q+\epsilon)(q^3+\epsilon)$&$q(q+\epsilon)(q^3+\epsilon)$\\
\hline
$\varphi_{1a},\varphi_{1b}$&$(q^2-1)(q^4+q^2+1)$&$(q^2-1)(q^4+q^2+1)$\\\hline
$\varphi_{2a},\varphi_{2b}$&$(q-1)^2(q^4+q^2+1)$&$(q-1)^2(q^4+q^2+1)$\\\hline
\end{tabular}
\end{center}
\footnotesize where:
\begin{enumerate}
  \item[(i)] $0\leq\alpha\leq q-1$ if $3\nmid q$ and $0\leq\alpha\leq 2q$ if $3\mid q$,
  \item[(ii)] $0\leq\beta\leq\frac{1}{3}(q+2)$,
  \item[(iii)] $1\leq\gamma\leq \frac{1}{3}(q+1)$.
\end{enumerate}
\normalsize Therefore, if $3\nmid q$ then $\varphi_{12}(1)\geq
\frac{1}{3}(q-1)^2(q+1)(q^3+2q^2+q+3)$ and if $3\mid q$ then
$\varphi_{12}(1)\geq\frac{1}{3}(q-1)^2(q^3+2q^2+4q+3)$. Moreover,
when $4\mid (q+1)$ and $q\equiv -1 (\bmod$ $3)$, we have
$\varphi_{31}(1)\geq \frac{1}{3}(q-1)^2(q^2+q+1)(2q^2+2q+3)$ (cf.
\cite{HS2}).

\medskip
\centerline {{\sc {Table III.}} The degrees of irreducible
$3$-Brauer characters of $G_2(q)$ with $3\nmid q$}

\begin{center}
\begin{tabular}{|c|c|c|} \hline
Character & \begin{tabular} {c}$3\mid (q-1)$, $0\leq\alpha\leq1$\\
$0\leq\beta\leq q-2$, $1\leq \gamma\leq q+1$\end{tabular} &
\begin{tabular} {c}$3\mid (q+1)$, $1\leq \alpha\leq q+1$\\$1\leq\beta\leq q-1$, $1\leq
\gamma\leq \frac{1}{2}q$\end{tabular}\\\hline
 $\varphi_{12}$& \begin{tabular} {c}$\frac{1}{6}(q-1)^2(6q^4+(9-\beta$\\$-2\gamma)
 q^3+(9+\beta-4\gamma)q^2$\\$+(9-\beta-2\gamma)q+6)$\end{tabular}&
 \begin{tabular} {c}$\frac{1}{6}(q-1)^2(6q^4+(11-\alpha-2\beta$\\$-3\gamma)
 q^3+(13+\alpha-4\beta-3\gamma)q^2$\\$+(11-\alpha-2\beta-3\gamma)q+6)$\end{tabular}\\ \hline
$\varphi_{14}$&\begin{tabular}{c}$\frac{1}{6}q(q^2-q+1)((1-\alpha)q^2$\\$
+(4+2\alpha)q+(1-\alpha))$\end{tabular}&
$\frac{1}{6}(q^2-1)(q^3+3q^2-q+6)$\\\hline
$\varphi_{15}$&$\frac{1}{2}(q^5+q^4+q^2+q-2)$&$\frac{1}{2}q(q+1)^2(q^2-q+1)$\\\hline
$\varphi_{16}$&$q^3$&$q^3-1$\\\hline
$\varphi_{21}$&$q^2(q^4+q^2+1)$&$(q-1)^2(q^4+q^2+1)$\\\hline
$\varphi_{22}$&$q^4+q^2+1$&$q^4+q^2+1$\\\hline $\varphi_{23},
\varphi_{24}$&$q(q^4+q^2+1)$&$(q-1)(q^4+q^2+1)$\\\hline
$\varphi_{1a},
\varphi_{1b}$&$q(q+1)(q^4+q^2+1)$&$(q^2-1)(q^4+q^2+1)$\\\hline
$\varphi_{2a},
\varphi_{2b}$&$q(q-1)(q^4+q^2+1)$&$(q-1)^2(q^4+q^2+1)$\\\hline
\end{tabular}
\end{center}
\footnotesize where:
\begin{enumerate}
  \item[(i)] when $3\mid (q-1)$, $\varphi_{14}(1)\in\{\frac{1}{6}q(q^2-q+1)(q^2+4q+1)$,
  $q^2(q^2-q+1)\}$ and \\$\varphi_{12}\geq \frac{1}{2}(q-1)^2(q^4+2q^3+3q+2)$,
  \item[(ii)] when $3\mid(q+1)$, $\varphi_{12}\geq \frac{1}{4}(q-1)^2(q+2)^2(q^2+q+1)$ (cf. \cite{HS1}).
\end{enumerate} \normalsize

By comparing the degrees of characters directly, we easily get the
value of $\mathfrak{d}_\ell(G)$. We include the case when $\ell\nmid
|G|$ in the following formula, where $q\geq 5$ and $\ell\nmid q$.
\begin{equation}\label{d1l}
\mathfrak{d}_\ell(G)= \left\{\begin {array}{ll}
q^3+1, & q\equiv 1 (\bmod$ $3)$, $\ell\neq 3,\\
q^3, & q\equiv 1 (\bmod$ $3)$, $\ell= 3,\\
q^3-1, & q\equiv 2 (\bmod$ $3)$, $\forall$ $\ell,\\
q^4+q^2, & q \equiv 0 (\bmod$ $3)$, $\ell=2,\\
q^4+q^2+1, & q\equiv 0 (\bmod$ $3)$, $\ell\neq2,
 \end {array} \right.
 \end{equation}

When $\ell\nmid |G|$, we know that $\IBR_\ell(G)=\Irr(G)$ and the
value of $\mathfrak{d}_{2,\mathbb{C}}(G)$ is given in formula
(\ref{second}). Formula (\ref{second}) and direct comparison of
degrees of $\ell$-Brauer characters when $\ell\mid |G|$ yields the
following theorem. We omit the details of this direct computation.

\begin{theorem} \label{d2l} Let $G=G_2(q)$ with $q\geq5$. We have
$$\mathfrak{d}_{2,2}(G)= \left\{\begin{array}{ll}
\frac{1}{6}q(q-1)^2(q^2-q+1), & p=3 \text{ or } q=5$, $7,\\
q^4+q^2, & p\geq5$, $q\geq11,
 \end{array} \right. $$
$$\mathfrak{d}_{2,3}(G) \left\{\begin {array}{ll}
=\frac{1}{6}q(q-1)^2(q^2-q+1), & q=5$, $7$, \text{or} $p=2$, $q\equiv -1(\bmod$ $3),\\
=q^4+q^2+1, & p\geq 5$, $q\equiv -1(\bmod$ $3)$\text{ and }$q\geq11, \\
\geq q^4-q^3+q^2, & q \geq13$\text{ and }$q\equiv1(\bmod$ $3),
 \end {array} \right. $$
and if $\ell=0$ or $\ell\geq5$, $\ell\nmid q$ then
$$\mathfrak{d}_{2,\ell}(G)=\mathfrak{d}_{2,\mathbb{C}}(G)= \left\{\begin {array}{ll}
\frac{1}{6}q(q-1)^2(q^2-q+1), & p=2$, $3$\text{ or }$q=5$, $7,\\
q^4+q^2+1, & p\geq5$, $q\geq11.
 \end {array} \right. $$
Moreover, $\IBR_\ell(G)$ contains a unique character denoted by
$\psi$ of degree larger than $1$ but less than
$\mathfrak{d}_{2,\ell}(G)$ and this character has degree
$\mathfrak{d}_{\ell}(G)$.
\end{theorem}

From the results about blocks and Brauer trees of $G_2(q)$ in
\cite{H1}, \cite{HS1}, \cite{HS2}, \cite{S1} and \cite{S2}, we see
that if $3\nmid q$ then $\psi=\widehat{X_{32}}$ in all cases except
when $\ell=3$ and $q\equiv1 (\bmod$ $3)$ where
$\psi=\widehat{X_{32}}-\widehat{{1}_G}$. If $3\mid q$ then
$\psi=\widehat{X_{22}}$ except when $\ell=2$ where
$\psi=\widehat{X_{22}}-\widehat{1_G}$. We also notice that
$\varphi_{18}=\widehat{X_{18}}$ is irreducible in all cases.

\section{Proof of main results}

In next lemmata, we use the notation in \cite{CR}, \cite{EY},
\cite{G}, \cite{H2}, and \cite{SF}.

\begin{lemma}\label{lemma}
\begin{enumerate}
  \item[(i)] When $q\equiv-1(\bmod$ $3)$, the character $X_{32}|_{SL_3(q)}$ is
irreducible and equal to $\chi_{q^3-1}^{(\frac{q^2-1}{3})}$.
  \item[(ii)] When $q\equiv1(\bmod$ $3)$, the character $X_{32}|_{SU_3(q)}$ is
irreducible and equal to $\chi_{q^3+1}^{(\frac{q^2-1}{3})}$.
\end{enumerate}
\end{lemma}

\begin{proof}
(i) First let us consider the case when $q$ is odd and
$q\equiv-1(\bmod$ $3)$. We need to find the fusion of conjugacy
classes of $SL_3(q)$  in $G$. By \cite[p. 487]{SF}, $SL_3(q)$ has
the following conjugacy classes: $C_1^{(0)}$, $C_2^{(0)}$,
$C_3^{(0,0)}$, $C_4^{(k)}$, $C_5^{(k)}$, $C_6^{(k,l,m)}$,
$C_7^{(k)}$ and $C_8^{(k)}$. We denote by $K(x)$ the conjugacy class
containing an element $x\in G$. For the representatives of the
conjugacy classes of G, we refer to \cite[p. 396-398]{CR}. Then we
have $C_1^{(0)}\subseteq K(1)$, $C_2^{(0)}\subseteq K(u_1)$,
$C_3^{(0,0)}\subseteq K(u_6)$, $C_4^{(k)}\subseteq K(k_2)\cup
K(h_{1b})$, $C_5^{(k)}\subseteq K(k_{2,1})\cup K(h_{1b,1})$,
$C_6^{(k,l,m)}\subseteq K(h_1)$, $C_7^{(k)}\subseteq K(h_b)\cup
K(h_{2b})$ and $C_8^{(k)}\subseteq K(h_3)$. Taking the values of
${X_{32}}$ and $\chi_{q^3-1}^{(\frac{q^2-1}{3})}$ on these classes,
we get
\begin{center}
\begin{tabular} {l}
$\chi_{q^3-1}^{(\frac{q^2-1}{3})}(C_1^{(0)})=X_{32}(1)=q^3-1$,\\
$\chi_{q^3-1}^{(\frac{q^2-1}{3})}(C_2^{(0)})=X_{32}(u_1)=-1$,\\
$\chi_{q^3-1}^{(\frac{q^2-1}{3})}(C_3^{(0,0)})=X_{32}(u_6)=-1$,\\
$\chi_{q^3-1}^{(\frac{q^2-1}{3})}(C_4^{(k)})=X_{32}(k_2)=X_{32}(h_{1b})=q-1$,\\
$\chi_{q^3-1}^{(\frac{q^2-1}{3})}(C_5^{(k)})=X_{32}(k_{2,1})=X_{32}(h_{1b,1})=-1$,\\
\end{tabular}

\begin{tabular} {l}
$\chi_{q^3-1}^{(\frac{q^2-1}{3})}(C_6^{(k,l,m)})=X_{32}(h_1)=0$,\\
$\chi_{q^3-1}^{(\frac{q^2-1}{3})}(C_7^{(k)})=X_{32}(h_b)=X_{32}(h_{2b})=-(\omega^k+w^{-k})$,\\
$\chi_{q^3-1}^{(\frac{q^2-1}{3})}(C_8^{(k)})=X_{32}(h_3)=0,$
\end{tabular}
\end{center}
where $\omega$ is a primitive cubic root of unity. Therefore,
${X_{32}}|_{SL_3(q)}= \chi_{q^3-1}^{(\frac{q^2-1}{3})}$. The case
$q$ even is proved similarly by using \cite{EY}.

(ii) This part is proved similarly by using \cite[p. 396-398]{CR}
and \cite[p. 565]{G}.
\end{proof}

\begin{lemma}\label{lemmacase3}
Suppose that $q\equiv-1(\bmod$ $3)$. Then
$\widehat{\chi_{q^3-1}^{(\frac{q^2-1}{3})}}\in \IBR_{\ell}(SL_3(q))$
for $\ell\nmid q$.
\end{lemma}

\begin{proof}
Since $q\equiv-1(\bmod$ $3)$, $GL_3(q)=SL_3(q)\times Z(GL_3(q))$
with $Z(GL_3(q))\simeq \mathbb{Z}_{q-1}$. Let $V$ be a
$\mathbb{C}GL_3(q)$-module affording the irreducible character
$\chi_{q^3-1}^{(\frac{q^2-1}{3})}\cdot{1}_{\mathbb{Z}_{q-1}}\in
\Irr(GL_3(q))$. We have $\dim V=q^3-1$. According to \cite[\S4]{GT},
$V$ has the form $(S_\mathbb{C}(s, (1))\circ
S_\mathbb{C}(t,(1)))\uparrow G$ where $s\in \mathbb{F}^\times_q$ and
$t$ has degree $2$ over $\mathbb{F}_q$. Using Corollary 2.7 of
\cite{GT}, since 1 and 2 are coprime, we see that $V$ is irreducible
in any cross characteristic. This implies that $V|_{SL_3(q)}$ is
also irreducible in any cross characteristic and therefore
$\widehat{\chi_{q^3-1}^{(\frac{q^2-1}{3})}}\in \IBR_\ell(SL_3(q))$
for every $\ell\nmid q$.
\end{proof}

Note: $2$-Brauer characters of $SU_3(q)$ were not considered in
\cite{G}.

\begin{lemma}\label{lemmacase4}
Suppose that $q\equiv1(\bmod$ $3)$ and $q$ is odd. Then
$\widehat{\chi_{q^3+1}^{(\frac{q^2-1}{3})}}\in\IBR_2(SU_3(q))$.
\end{lemma}

\begin{proof}
We denote the character $\chi_{q^3+1}^{(\frac{q^2-1}{3})}$ by $\rho$
for short. Since $q\equiv1(\bmod$ $3)$, $GU_3(q)=SU_3(q)\times
Z(GU_3(q))$ with $Z(GU_3(q))\simeq Z_{q+1}$. Hence, $SU_3(q)$ has
the same degrees of irreducible Brauer characters as $GU_3(q)$ does.
It is shown in \cite{TZ} that every $\varphi\in \IBR_2(GU_3(q))$
either lifts to characteristic 0 or $\varphi(1)=q(q^2-q+1)-1$. This
and the character table of $SU_3(q)$ (cf. \cite{G}) gives the
possible values for degrees of $2$-Brauer characters of $SU_3(q)$:
$1$, $q^2-q$, $q^2-q+1$, $q(q^2-q+1)-1$, $q(q^2-q+1)$,
$(q-1)(q^2-q+1)$, $q^3$, $q^3+1$, and $(q+1)^2(q-1)$.

Assume $\widehat{\rho}$ is reducible. Then $\widehat{\rho}$ is the
sum of more than one irreducible characters of $SU_3(q)$. Hence,
$\rho(1)=q^3+1$ is the sum of more than one of the values listed
above. It follows that $\widehat{\rho}$ must include a character of
degree either $1$, or $q^2-q$, or $q^2-q+1$. Once again, this
character lifts to a complex character that we denote by $\alpha$.
Clearly, $\rho$ and $\alpha$ belong to the same $2$-block of
$SU_3(q)$. We will use central characters to show that this can not
happen.

Let $R$ be the full ring of algebraic integers in $\mathbb{C}$ and
$\pi$ a maximal ideal of $R$ containing $2R$. It is known that
$\alpha$ and $\rho$ are in the same 2-block if and only if
\begin{equation}\label{pi} \omega_\rho(K)-\omega_\alpha(K)\in \pi
\end{equation}
where $K$ is any class sum and $\omega_\chi$ is the central
character associated with $\chi$. The value of $\omega_\chi$ on a
class sum is
$$\omega_\chi(K)=\frac{\chi(g)|\mathcal{K}|}{\chi(1)}$$ where $\mathcal{K}$
is the conjugacy class with class sum $K$ and $g$ is an element in
$\mathcal{K}$. Therefore, (\ref{pi}) implies that
\begin{equation}\label{pi1}
\frac{\rho(g)}{\rho(1)}|g^G|-\frac{\alpha(g)}{\alpha(1)}|g^G|\in \pi
\end{equation}
where $|g^G|$ denotes the length of the conjugacy class containing
$g\in G$.

First assume that $\alpha$ is the trivial character. In (\ref{pi1}),
take $g$ to be any element in the conjugacy class $C_7^{(1)}$ (cf.
\cite{G}), we have
$\frac{\rho(g)}{\rho(1)}|g^G|-\frac{\alpha(g)}{\alpha(1)}|g^G|=
\frac{-1}{q^3+1}q^3(q^3+1)-\frac{1}{1}q^3(q^3+1)=-q^3-q^3(q^3+1)\in
\pi$. Note that $\pi\cap \mathbb{Z}=2\mathbb{Z}$. Since
$-q^3-q^3(q^3+1)$ is an odd number, we get a contradiction.

Secondly, if $\alpha(1)=q^2-q$ then $\alpha=\chi_{q^2-q}$ as denoted
in \cite{G}. In (\ref{pi1}), take $g$ to be any element in the
conjugacy class $C_7^{(1)}$, we have
$\frac{\rho(g)}{\rho(1)}|g^G|-\frac{\alpha(g)}{\alpha(1)}|g^G|=
\frac{-1}{q^3+1}q^3(q^3+1)-\frac{0}{q^2-q}q^3(q^3+1)=-q^3\in \pi$.
Since $-q^3$ is an odd number, we again get a contradiction.

Finally, if $\alpha(1)=q^2-q+1$ then $\alpha=\chi_{q^2-q+1}^{(u)}$
for some $1\leq u\leq q$. In (\ref{pi1}), take $g$ to be any element
in the conjugacy class $C_7^{(q+1)}$. Then we have
$\frac{\rho(g)}{\rho(1)}|g^G|-\frac{\alpha(g)}{\alpha(1)}|g^G|=
\frac{-1}{q^3+1}q^3(q^3+1)-\frac{1}{q^2-q+1}q^3(q^3+1)=-q^3-q^3(q+1)$,
which is an odd number, a contradiction.
\end{proof}

\begin{theorem} [Reduction Theorem]
\label{red} Let $G=G_2(q),$ $q=p^n$, $q\geq 5$, and $p$ a prime
number. Let $\Phi$ be an absolutely irreducible representation of
$G$ in cross characteristic $\ell$ and $M$ a maximal subgroup of
$G$. Assume that $\deg(\Phi)>1$ and $\Phi|_ M$ is also absolutely
irreducible. Then $M$ is $G$-conjugate to one of the following
groups:
\begin{enumerate}
  \item[(i)] maximal parabolic subgroups $P_a, P_b$,
  \item[(ii)] $SL_3(q):2$, $SU_3(q):2$,
  \item[(iii)] $G_2(q_0)$ with $q=q_0^2$, $p\neq 3$.
\end{enumerate}
\end{theorem}

\begin{proof}
By Lemmata \ref{1} and \ref{2}, we have $\mathfrak{d}_{\ell}(G)\leq
\deg(\Phi)\leq \mathfrak{m}_\ell(M)\leq
\mathfrak{m}_{\mathbb{C}}(M)\leq \sqrt{|M|}$. Hence,
$\mathfrak{d}_\ell(G)\leq\sqrt{|M|}$. Moreover, from formulas
(\ref{d1c}) and (\ref{d1l}), we have $\mathfrak{d}_\ell(G)\geq
q^3-1$ if $3\nmid q$ and $\mathfrak{d}_\ell(G)\geq q^4+q^2$ if
$3\mid q$ for every $\ell\nmid q$. Therefore,
\begin{equation}\label{root}
\sqrt{|M|}\geq \left\{\begin {array}{ll} q^3-1 &\text{if } 3\nmid q,\\
q^4+q^2 &\text{if } 3\mid q.
\end {array} \right.
\end{equation}

Here, we will only give the proof for the case $p\geq 5$. The proofs
for $p=2$ and $p=3$ are similar. According to \cite{K1}, if $M$ is a
maximal subgroup of  $G$ with $p\geq5$ then $M$ is $G$-conjugate to
one of the following groups:
\begin{enumerate}
  \item[1)] $P_a, P_b$, maximal parabolic subgroups,
  \item[2)] $(SL_2(q)\circ SL_2(q))\cdot 2$, involution centralizer,
  \item[3)] $2^3 \cdot L_3(2)$, only when $p=q$,
  \item[4)] $SL_3(q):2$, $SU_3(q):2$,
  \item[5)] $G_2(q_0)$, $q=q_o^\alpha$, $\alpha$ prime,
  \item[6)] $PGL_2(q)$, $p\geq 7$, $q\geq 11$,
  \item[7)] $L_2(8)$, $p\geq 5$,
  \item[8)] $L_2(13)$, $p \neq 13$,
  \item[9)] $G_2(2)$, $q=p\geq 5$,
  \item[10)] $J_1$, $q=11$.
\end{enumerate}
Consider for instance the case 5) with $\alpha\geq 3$. Then
$\sqrt{|M|}=\sqrt{q_0^6(q_0^6-1)(q_0^2-1)}<q_0^7\leq q^{7/3}<q^3-1$
for every $q\geq5$. This contradicts (\ref{root}). The cases 2), 3),
6)-10) are excluded similarly.
\end{proof}

{\bf Proof of Theorem 1.1.} By the Reduction Theorem, $M$ is
$G$-conjugate to one of the following subgroups:
\begin{enumerate}
  \item[(i)] maximal parabolic subgroups $P_a, P_b$,
  \item[(ii)] $SL_3(q):2$, $SU_3(q):2$,
  \item[(iii)] $G_2(q_0)$ with $q=q_0^2$, $p\neq 3$.
\end{enumerate}

Now we will proceed case by case.

\medskip

$\mathbf{Case}$ $\mathbf{1}$: $M=P_a$.

Let $Z:=Z(P_a')$, the center of the derived subgroup of $P_a$. We
know that $P_a$ is the normalizer of $Z$ in $G$ and therefore $Z$ is
nontrivial. In fact, $Z$ is a long-root subgroup of $G$. Let $V$ be
an irreducible $G$-module affording the character $\varphi$. By
Lemma \ref{lemmaforparabolica1}, $Z$ must have nonzero fixed points
on $V$. In other words, $C_V(Z)=\{v\in V\mid a(v)=v$ for every $a\in
Z\}\neq0$. Therefore $V|_{P_a}$ is reducible by Lemma \ref{lemma for
Parabolic}(i).

\medskip

$\mathbf{Case}$ $\mathbf{2}$: $M=P_b$.

Using the results about character tables of $P_b$ in \cite{AH},
\cite{E}, and \cite{EY}, we have $\mathfrak{m}_{\mathbb{C}}(P_b)=
q(q-1)(q^2-1)$ for $q\geq 5$. Therefore, if $\varphi|_{P_b}$ is
irreducible then $\varphi(1)\leq q(q-1)(q^2-1)$. If $3\mid q$ then
$\mathfrak{d}_{\ell}(G)\geq q^4+q^2$ because of formula (\ref{d1l}).
Then we have $\mathfrak{d}_{\ell}(G)\geq q^4+q^2> q(q-1)(q^2-1)\geq
\varphi(1)$ and this cannot happen. So $q$ must be coprime to $3$.

It is easy to check that $\mathfrak{m}_{\mathbb{C}}(P_b)<
\mathfrak{d}_{2,\ell}(G)$ for every $q\geq8$. Therefore, when
$q\geq8$, the inequality $\varphi(1)\leq q(q-1)(q^2-1)$ can hold
only if $\varphi$ is the nontrivial character of smallest degree.
When $q=5$ or $7$, checking directly, we see that besides the
nontrivial character of smallest degree, $\varphi$ can be
$\varphi_{18}=\widehat{X_{18}}$ of degree
$\frac{1}{6}q(q-1)^2(q^2-q+1)$. Recall that $|P_b|=q^6(q^2-1)(q-1)$,
which is not divisible by $X_{18}(1)$. Therefore ${X_{18}|_{P_b}}$
is reducible and so is $\varphi_{18}|_{P_b}$ for $q=5$ or $7$. We
have shown that the unique possibility for $\varphi$ is the
character $\psi$ of smallest degree when $3\nmid q$.

Now note that $X_{32}(1)=(q^3+\epsilon)\nmid |P_b|$ and therefore
$X_{32}|_{P_b}$ must be reducible. If $\ell=3$ and $q\equiv1 (\bmod$
$3)$ then $\psi=\widehat{X_{32}}-\widehat{1_G}$. Assume that
$\psi|_{P_b}=\widehat{X_{32}}|_{P_b}-\widehat{1_{P_b}}$ is
irreducible. The reducibility of $X_{32}|_{P_b}$ and the
irreducibility of $\widehat{X_{32}}|_{P_b}-\widehat{1_{P_b}}$ imply
that ${X_{32}}|_{P_b}=\lambda+\mu$ where
$\widehat{\lambda}=\widehat{1_{P_b}}$, $\mu\in \Irr(P_b)$ and
$\widehat{\mu}\in \IBR_3(P_b)$. We then have
$\mu(1)=X_{32}(1)-1=q^3$. Inspecting the character tables of $P_b$
given in \cite{AH} and \cite{EY}, we see that there is no
irreducible complex character of $P_{b}$ of degree $q^3$ and we get
a contradiction. If $\ell\neq3$ or $q$ is not congruent to 1 modulo
3 then $\psi=\widehat{X_{32}}$. Therefore
$\psi|_{P_b}=\widehat{X_{32}}|_{P_b}$, which is reducible as noted
above.

\medskip
$\mathbf{Case}$ $\mathbf{3}$: $M=SL_3(q):2$.

From \cite{SF}, we know that
$\mathfrak{m}_{\mathbb{C}}(SL_3(q))=(q+1)(q^2+q+1)$ for every $q\geq
5$. Therefore,
$\mathfrak{m}_{\mathbb{C}}(SL_3(q):2)\leq2(q+1)(q^2+q+1)$. Since
$\varphi|_{SL_3(q):2}$ is irreducible, $\varphi(1)\leq
2(q+1)(q^2+q+1)$. Similarly as the previous case, we have
$q^4+q^2>2(q+1)(q^2+q+1)$ for every $q\geq5$ and so $q$ is not
divisible by $3$.

By Theorem \ref{d2l}, the inequality $\varphi(1)\leq
2(q+1)(q^2+q+1)$ can hold only if $\varphi$ is the character $\psi$
of smallest degree or $\varphi_{18}=\widehat{X_{18}}$ when $q=5$. We
have $X_{18}(1)=280$, $|SL_3(5):2|=744,000$ and therefore
$X_{18}(1)\nmid |SL_3(q):2|$ when $q=5$. Hence,
${X_{18}|_{SL_3(q):2}}$ is reducible and so is
$\varphi_{18}|_{SL_3(q):2}$ when $q=5$. Again, the unique
possibility for $\varphi$ is $\psi$ when $3\nmid q$.

If $\ell=3$ and $q\equiv 1(\bmod$ $3)$ then
$\psi=\widehat{X_{32}}-\widehat{1_G}$. Assume that
$\psi|_{SL_3(q):2}$ is irreducible. Let $V$ be an irreducible
$G$-module in characteristic $3$ affording the character $\psi$.
Then $V|_{SL_3(q):2}$ is an irreducible $(SL_3(q):2)$-module. Let
$\sigma$ be a generator for the multiplicative group
$\mathbb{F}_q^\times$ and $I$ be the identity matrix in $SL(3,
\mathbb{F}_q)$. Consider the matrix $T=\sigma^{\frac{q-1}{3}}\cdot
I$. We have $<T>= Z(SL_3(q))$ and hence $<T>~\trianglelefteq
SL_3(q):2$. Since $\ord(T)=3$ and $<T>~\leqslant O_3(SL_3(q):2)$, it
follows that $O_3(SL_3(q):2)$ is nontrivial. By Lemma \ref{lemma for
Parabolic}(ii), $V|_{SL_3(q):2}$ is reducible and we get a
contradiction. If $\ell\neq 3$ and $q\equiv1 (\bmod$ $3)$ then
$\psi=\widehat{X_{32}}$. Note that $X_{32}(1)=q^3+1$ and
$|SL_3(q):2|=2q^3(q^3-1)(q^2-1)$. Therefore $X_{32}(1)\nmid
|SL_3(q):2|$ for every $q\geq5$. It follows that
$X_{32}|_{SL_3(q):2}$ as well as $\psi|_{SL_3(q):2}$ are reducible.

It remains to consider $q\equiv -1 (\bmod$ $3)$ and therefore
$\psi=\widehat{X_{32}}$. By Lemmata \ref{lemma} and
\ref{lemmacase3}, we get that
$\psi|_{SL_3(q)}=\widehat{X_{32}}|_{SL_3(q)}
=\widehat{\chi_{q^3-1}^{(\frac{q^2-1}{3})}}\in\IBR_{\ell}(SL_3(q))$
for $\ell\nmid q$, as we claim in the item (i) of Theorem \ref{main
theorem}.

\medskip
$\mathbf{Case}$ $\mathbf{4}$: $M=SU_3(q):2$.

According to \cite{SF},
$\mathfrak{m}_{\mathbb{C}}(SU_3(q))=(q+1)^2(q-1)$ and therefore
$\mathfrak{m}_{\mathbb{C}}(SU_3(q):2)\leq 2(q+1)^2(q-1)$ for every
$q\geq 5$. Hence $\varphi(1)\leq 2(q+1)^2(q-1)$ by the
irreducibility of $\varphi|_{SU_3(q):2}$. Again, $q$ must be coprime
to $3$.

By Theorem \ref{d2l}, the inequality $\varphi(1)\leq 2(q+1)^2(q-1)$
can hold only if $\varphi$ is the character $\psi$ of smallest
degree or $\varphi_{18}$ of degree $\frac{1}{6}q(q-1)^2(q^2-q+1)$
when $q=5$. By \cite{G}, the degrees of irreducible characters of
$SU_3(5)$ are: 1, 20, 125, 21, 105, 84, 126, 144, 28 and 48. When
$q=5$, $X_{18}(1)=280$. Therefore, $X_{18}|_{SU_3(5)}$ is the sum of
at least $3$ irreducible characters. Since $SU_3(5)$ is a normal
subgroup of index $2$ of $SU_3(5):2$, by Clifford's theorem,
$X_{18}|_{SU_3(5):2}$ is reducible when $q=5$. This implies that
$\varphi_{18}$ is also reducible when $q=5$. So, the unique
possibility for $\varphi$ is $\psi$ when $3\nmid q$.

If $q\equiv -1 (\bmod$ $3)$ then $\psi=\widehat{X_{32}}$. Note that
$|SU_3(q):2|=2q^3(q^3+1)(q^2-1)$, which is not divisible by $q^3-1$
for every $q\geq 5$. Hence $X_{32}|_{SU_3(q):2}$ is reducible and so
is $\psi|_{SU_3(q):2}$. It remains to consider $q\equiv 1 (\bmod$
$3)$. First, if $\ell=2$ then $\psi=\widehat{X_{32}}$. By Lemmata
\ref{lemma} and \ref{lemmacase4}, we have
$\psi|_{SU_3(q)}=\widehat{X_{32}}|_{SU_3(q)}=\widehat{\chi_{q^3+1}^{(\frac{q^2-1}{3})}}\in
\IBR_2(SU_3(q))$. Next, if $\ell= 3$ then
$\psi=\widehat{X_{32}}-\widehat{1_G}$. By Lemma \ref{lemma}, we have
$\psi|_{SU_3(q)}=\widehat{X_{32}}|_{SU_3(q)}-\widehat{1}_{SU_3(q)}=
\widehat{\chi_{q^3+1}^{(\frac{q^2-1}{3})}}-\widehat{1}_{SU_3(q)}=
\widehat{\chi_{q^3}}$ which is an irreducible 3-Brauer character of
$SU_3(q)$ (cf. \cite[p. 573]{G}). Finally, if $\ell\neq 2,3$ then
$\psi=\widehat{X_{32}}$. By Lemma \ref{lemma}, we have
$\widehat{X_{32}}|_{SU_3(q)}=\widehat{\chi_{q^3+1}^{(\frac{q^2-1}{3})}}$
which is irreducible by \cite{G}. We have shown that when
$q\equiv1(\bmod$ $3)$, $\psi|_{SU_3(q):2}$ is irreducible, as we
claim in the item (ii) of Theorem \ref{main theorem}.

\medskip
$\mathbf{Case}$ $\mathbf{5}$: $M=G_2(q_0)$ with $q=q_0^2$, $3\nmid
q$.

Since $\varphi|_M$ is irreducible, $\varphi(1)\leq
\sqrt{|G_2(q_0)|}=\sqrt{q_0^6(q_0^6-1)(q_0^2-1)}<\sqrt{q^7}<\mathfrak{d}_{2,\ell}(G)$
for every $q\geq5$. Therefore, by Theorem \ref{d2l}, the unique
possibility for $\varphi$ is $\psi$. Since $q=q_0^2$, $q\equiv 1
(\bmod$ $3)$.

Recall that $X_{32}(1)=q^3+1$ and
$|G_2(q_0)|=q_0^6(q_0^6-1)(q_0^2-1)=q^3(q^3-1)(q-1)$. It is easy to
see that $(q^3+1)\nmid q^3(q^3-1)(q-1)$ for every $q\geq 5$. So
$X_{32}|_{G_2(q_0)}$ is reducible. First we consider the case
$\ell=3$. Then $\psi=\widehat{X_{32}}-\widehat{1_G}$. Assume that
$\psi|_{G_2(q_0)}=\widehat{X_{32}}|_{G_2(q_0)}-\widehat{1}_{G_2(q_0)}$
is irreducible. Then ${X_{32}}|_{G_2(q_0)}=\lambda+\mu$ where
$\widehat{\lambda}=\widehat{1}_{G_2(q_0)}$, $\mu\in \Irr(G_2(q_0))$
and $\widehat{\mu}\in \IBR_3(G_2(q_0))$. We then have
${\mu}(1)=X_{32}(1)-1=q^3=q_0^6$. So $\mu$ is the Steinberg
character of $G_2(q_0)$. From \cite{HS1}, we know that the reduction
modulo $3$ of the Steinberg character is reducible, a contradiction.
Now we can assume $\ell\neq3$ and therefore $\psi=\widehat{X_{32}}$.
From the reducibility of $X_{32}|_{G_2(q_0)}$ as noted above,
$\psi|_{G_2(q_0)}$ is reducible. \hfill $\Box$

\section{Small groups}

In this section, we mainly use results and notation of \cite{Atlas1}
and \cite{Atlas2}.

\begin{lemma}\label{smalllemma1}
Theorem \ref{main theorem2} holds in the case $G=G_2(3)$,
$M=U_3(3):2$.
\end{lemma}

\begin{proof}
According to \cite[p. 14]{Atlas1}, we have
$\mathfrak{m}_{\mathbb{C}}(U_3(3))=32$ and
$\mathfrak{m}_{\mathbb{C}}(U_3(3):2)=64$. Thus, if $\varphi|_M$ is
irreducible then $\varphi(1)\leq 64$. Inspecting the character
tables of $G_2(3)$ in \cite[p. 60]{Atlas1} and \cite[p. 140, 142,
143]{Atlas1}, we see that $\varphi(1)=14$ or $64$.

Note that $G_2(3)$ has a unique irreducible complex character of
degree $14$ which is denoted by $\chi_2$ and every reduction modulo
$\ell\neq 3$ of $\chi_2$ is still irreducible. Now we will show that
$\chi_2|_{U_3(3)}=\chi_6$, which is the unique irreducible character
of degree $14$ of $U_3(3)$. Suppose that
$\chi_2|_{U_3(3)}\neq\chi_6$, then $\chi_2|_{U_3(3)}$ is reducible
and it is the sum of more than one irreducible characters of degree
less than $14$. Note that $U_3(3)$ has exactly one conjugacy class
of elements of order $6$, which is denoted by $6A$. If
$\chi_2|_{U_3(3)}$ is the sum of two irreducible characters, then
the degrees of these characters are both equal to $7$. So
$\chi_2|_{U_3(3)}(6A)$ is $0$, $2$ or $4$. This cannot happen since
the value of $\chi_2$ on any class of elements of order $6$ of
$G_2(3)$ is $1$ or $-2$. If $\chi_2|_{U_3(3)}$ is the sum of more
than two irreducible characters, then $\chi_2|_{U_3(3)}(6A)\geq 2$
which cannot happen, neither. So we have $\chi_2|_{U_3(3)}=\chi_6$.
We also see that every reduction modulo $\ell\neq 3$ of $\chi_6$ is
still irreducible. Hence, if $\ell\neq 3$ and $\varphi$ is the
irreducible $\ell$-Brauer character of $G_2(3)$ of degree $14$, then
$\varphi|_{U_3(3)}$ is irreducible and so is $\varphi|_{U_3(3):2}$.

By \cite[p. 14]{Atlas1}, $U_3(3):2$ has one irreducible complex
character of degree $64$ which we denote by $\chi$. Also, $G_2(3)$
has two irreducible complex characters of degree $64$ that are
$\chi_3$ and $\chi_4$ as denoted in \cite[p. 60]{Atlas1}. We will
show that $\chi_{3,4}|_{U_3(3):2}=\chi$. Note that the two conjugacy
classes of $(U_3(3):2)$-subgroups of $G_2(3)$ are fused under an
outer automorphism $\tau$ of $G_2(3)$, which stabilizes each of
$\chi_3$ and $\chi_4$. Hence, without loss we may assume that
$U_3(3):2$ is the one considered in \cite[p. 237]{E}. Checking
directly, it is easy to see that the values of $\chi_3$ and $\chi_4$
coincide with those of $\chi$ at every conjugacy classes except at
the classes of elements of order $3$. $U_3(3):2$ has two classes of
elements of order $3$, $3A$ and $3B$. Using \cite[p. 237]{E} to find
the fusion of conjugacy classes of $U_3(3)$ in $G_2(3)$ and the
values of $\chi_3$ and $\chi_4$ (which are $\theta_{12}(k)$ in
\cite{E}), we see that the classes $3A$, $3B$ of $U_3(3)$ are
contained in the classes $3A$, $3E$ of $G_2(3)$, respectively.
Moreover, $\chi_3(3A)=\chi_4(3A)=\chi(3A)=-8$ and
$\chi_3(3E)=\chi_4(3E)=\chi(3B)=-2$. We have shown that
$\chi_3|_{U_3(3):2}=\chi_4|_{U_3(3):2}=\chi$. By \cite[p. 140, 142,
143]{Atlas2}, the reductions modulo $\ell\neq3$ of $\chi_3$ as well
as $\chi_4$ are irreducible. Also, the reduction modulo $\ell$ of
$\chi$ is irreducible for every $\ell\neq 3$, $7$. Therefore,
$\widehat{\chi_3}|_{U_3(3):2}$ and $\widehat{\chi_4}|_{U_3(3):2}$
are equal and irreducible for every $\ell\neq 3$, $7$. When
$\ell=7$, both $\widehat{\chi_3}|_{U_3(3):2}$ and
$\widehat{\chi_4}|_{U_3(3):2}$ are reducible since
$\mathfrak{m}_7(U_3(3):2)=56$ and
$\widehat{\chi_3}(1)=\widehat{\chi_4}(1)=64$.
\end{proof}

\begin{lemma}\label{smalllemma2}
Theorem \ref{main theorem2} holds in the case $G=G_2(3)$,
$M=2^3\cdot L_3(2)$.
\end{lemma}

\begin{proof}
We have $|2^3\cdot L_3(2)|=1344$. So
$\mathfrak{m}_{\mathbb{C}}(2^3\cdot L_3(2))\leq \sqrt{1344}<37$ and
therefore the unique possibility for $\varphi$ is the reduction
modulo $\ell\neq3$ of the character $\chi_2$ of degree $14$.

When $\ell\nmid |G_2(3)|$, i.e. $\ell\neq 2$, $3$, $7$, and $13$,
from \cite{KT}, we know that $\widehat{\chi_2}|_{2^3\cdot L_3(2)}$
is irreducible. When $\ell=2$, we have $\mathfrak{m}_2(2^3\cdot
L_3(2))=\mathfrak{m}_2(L_3(2))\leq\sqrt{168}<13$. So
$\widehat{\chi_2}|_{2^3\cdot L_3(2)}$ is reducible when $\ell=2$.
When $\ell=13$, since $13\nmid |2^3\cdot L_3(2)|$,
$\widehat{\chi_2}|_{2^3\cdot L_3(2)}$ is still irreducible. The last
case we need to consider is $\ell=7$. Let $E=2^3\unlhd 2^3\cdot
L_3(2)$ which is an elementary abelian group of order $2^3$. Since
$\chi_2|_{E\cdot L_3(2)}$ is irreducible and $L_3(2)$ acts
transitively on $E\backslash \{1\}$ and $\Irr(E)\backslash
\{\mathbf{1}_E\}$,
$\chi_2|_{E}=2\cdot\sum_{\alpha\in\Irr(E)\backslash\{\mathbf{1}_E\}}
\alpha$. Let $I$ be the inertia group of $\alpha$ in $E\cdot
L_3(2)$. By Clifford's theorem, we have $\chi_2|_{E\cdot
L_3(2)}=\Ind_I^{E\cdot L_3(2)}(\rho)$ for some $\rho\in\Irr(I)$ and
$\rho|_E=2\alpha$. We also have $|I|=\frac{|E\cdot
L_3(2)|}{|\Irr(E)\backslash\{\mathbf{1}_E\}|}=2^6.3$, which is not
divisible by $7$. Thus, the reduction modulo $7$ of $\rho$ is
irreducible. Hence $\widehat{\chi_2}|_{E\cdot L_3(2)}$ is also
irreducible when $\ell=7$.
\end{proof}

\begin{lemma}\label{smalllemma3}
Theorem \ref{main theorem2} holds in the case $G=3\cdot G_2(3)$,
$M=3.P$ or $3.Q$ where $P$, $Q$ are maximal parabolic subgroups of
$G_2(3)$.
\end{lemma}

\begin{proof}
First, we consider $M=3.P$ where $P$ is one of two maximal parabolic
subgroups which is specified in \cite[p. 217]{E}. If
$\varphi|_{3.P}$ is irreducible then $\varphi(1)\leq
\mathfrak{m}_\mathbb{C}(3.P)\leq \sqrt{|(3.P)/Z(3.P)|}\leq
\sqrt{|P|}=\sqrt{11664}=108$. Inspecting the character tables (both
complex and Brauer) of $3\cdot G_2(3)$ in \cite{Atlas1} and
\cite{Atlas2} (note that we only consider faithful characters), we
have $\varphi(1)=27$ and $\varphi$ is actually the reduction modulo
$\ell\neq3$ of any of the four irreducible complex characters of
degree $27$ of $3\cdot G_2(3)$. From now on, we denote these
characters by $\chi_{24}$, $\overline{\chi_{24}}$ (corresponding to
the line $\chi_{24}$ in \cite[p. 60]{Atlas1}) and $\chi_{25}$,
$\overline{\chi_{25}}$ (corresponding to the line $\chi_{25}$ in
\cite[p. 60]{Atlas1}).

Now we will show that $\chi_{24}|_{3.P}$ is irreducible. Note that
if $g_1$ and $g_2$ are the pre-images of an element $g\in G_2(3)$
under the natural projection $\pi: 3\cdot G_2(3)\rightarrow G_2(3)$,
then $\chi_{24}(g_1)=\omega\chi_{24}(g_2)$ where $\omega$ is a cubic
root of unity. Therefore we have $[\chi_{24}|_{3.P},
\chi_{24}|_{3.P}]_{3.P}=\frac{1}{3\cdot |P|}\sum_{x\in 3.P}
\chi_{24}(x)\overline{\chi_{24}(x)}=\frac{1}{|P|}\sum_{g\in
P}\chi_{24}(\overline{g})\overline{\chi_{24}(\overline{g})}$, where
$\overline{g}$ is a fixed pre-image for each $g$. We choose
$\overline{g}$ so that the value of $\chi_{24}$ at $\overline{g}$ is
printed in \cite[p. 60]{Atlas1}. In \cite[p. 217]{E}, we have the
fusion of conjugacy classes of $P$ in $G_2(3)$. By comparing the
orders of centralizers of conjugacy classes of $G_2(3)$ in \cite[p.
239]{E} with those in \cite[p. 60]{Atlas1}, we can find a
correspondence between conjugacy classes of $G_2(3)$ in these two
papers. The length of each conjugacy class of $P$ can be computed
from \cite[p. 217, 218]{E}. All the above information is collected
in Table IV. From this table, we see that the value of $\chi_{24}$
is zero at any element $\overline{g}$ for which the order of $g$ is
$3$, $6$ or $9$. We denote by $\chi_{24}(\overline{X})$ the value
$\chi_{24}(\overline{g})$ for some $g\in X$, where $X$ is a
conjugacy class of $P$. Then we have
\begin{equation}\label{equation1}
\begin{tabular}{lll} $\sum_{g\in
P}\chi_{24}(\overline{g})\overline{\chi_{24}(\overline{g})}$&$=$&$|A_1|(\chi_{24}
(\overline{A_1}))^2+|B_{11}|(\chi_{24}(\overline{B_{11}}))^2$\\
&&$+|B_{21}|(\chi_{24}(\overline{B_{21}}))^2+|D_1|(\chi_{24}(\overline{D_1}))^2$\\
&&$+|D_2|(\chi_{24}(\overline{D_2}))^2+2|E_2(1)|(\chi_{24}(\overline{E_2(1)}))^2.$
\end{tabular}
\end{equation}
By \cite[p. 217]{E}, the class $D_1$ of $P$ is contained in the
class $D_{11}$ of $G_2(3)$. This class is the class $4A$ or $4B$
according to the notation in \cite[p. 60]{Atlas1}. First we assume
$D_{11}$ is $4A$. Then by looking at the values of one character of
degree $273$ of $G_2(3)$ both in \cite{E} and \cite{Atlas1}, it is
easy to see that $D_2\subset D_{12}=12A$ and $E_2(i)\subset E_2=8A$.
Therefore, formula (\ref{equation1}) becomes
\begin{center}
\begin{tabular}{lll}$\sum_{g\in
P}\chi_{24}(\overline{g})\overline{\chi_{24}(\overline{g})}$&$=$&
$27^2+324\cdot3^2+81\cdot3^2+486\cdot(-1)^2$\\
&&$+972\cdot2^2+2\cdot1458\cdot(-1)^2$\\
&$=$&$11,664$.
\end{tabular}
\end{center}
Next, we assume $D_{11}$ is $4B$. Then $D_{12}=12B$, $E_2=8B$ and
formula (\ref{equation1}) becomes
\begin{center}
\begin{tabular}{lll}$\sum_{g\in
P}\chi_{24}(\overline{g})\overline{\chi_{24}(\overline{g})}$&$=$&
$27^2+324\cdot3^2+81\cdot3^2+486\cdot 3^2+972\cdot0^2+2\cdot1458\cdot 1^2$\\
&$=$&$11,664$.
\end{tabular}
\end{center}
So in any case, we have $[\chi_{24}|_{3.P}, \chi_{24}|_{3.P}]_{3.P}=
\frac{1}{|P|}\sum_{g\in
P}\chi_{24}(\overline{g})\overline{\chi_{24}(\overline{g})}=1$. That
means $\chi_{24}|_{3.P}$ is irreducible.

\medskip
\centerline{{\sc {Table IV.}} Fusion of conjugacy classes of $P$ in
$G_2(3)$}

\begin{center}
\begin{tabular}{|c|c|c|c|} \hline
Fusion in \cite{E}& Corresponding class in \cite{Atlas1}&Length&
 Values of $\chi_{24}$\\\hline
 $A_1\subset A_1$&$1A$& $1$&$27$\\\hline
 $A_2\subset A_2$&&&\\
 $A_3\subset A_2$&&&\\
 $A_{41}\subset A_{31}$&&&\\
 $A_{42}\subset A_{32}$&&&\\
 $A_{43}\subset A_{32}$&&&\\
 $A_{51}\subset A_{41}$&&&\\
 $A_{52}\subset A_{42}$&$3A, 3B, 3C, 3D, 3E$&&$0$\\
 $A_{61}\subset A_{31}$&&&\\
 $A_{62}\subset A_{41}$&&&\\
 $A_{63}\subset A_{42}$&&&\\
 $A_{64}\subset A_{41}$&&&\\
 $A_{65}(t)\subset A_{32}, A_{41}, A_{42}$&&&\\
 $A_{66}(t)\subset A_{41}, A_{42}$&&&\\\hline
 $A_{71}\subset A_{51}$&&&\\
 $A_{72}\subset A_{52}$&$9A, 9B, 9C$&&$0$\\
 $A_{73}\subset A_{53}$&&&\\\hline
 $B_{11}\subset B_1$&$2A$&$324$&$3$\\ \hline
 \begin{tabular}{c}$B_{12}\subset B_2$\\
 $B_{13}\subset B_3$\\
 $B_{14}\subset B_4$\\
 $B_{15}\subset B_5$\\\end{tabular}
 &$6A, 6B, 6C, 6D$&&$0$\\\hline
 $B_{21}\subset B_1$&$2A$&$81$&$3$\\\hline
 \begin{tabular}{c}$B_{22}\subset B_2$\\
 $B_{23}\subset B_3$\\
 $B_{24}\subset B_4$\\
 $B_{25}\subset B_5$\\\end{tabular}&$6A, 6B, 6C, 6D$&&$0$\\\hline
 $D_{1}\subset D_{11}$&$4A, 4B$&$486$&$-1$, $3$\\\hline
 $D_{2}\subset D_{12}$&$12A, 12B$&$972$&$2$, $0$\\\hline
 \begin{tabular}{c}$E_{2}(i)\subset E_{2}$\\
 (two classes)\end{tabular}
 &$8A, 8B$&$1458$&$-1$, $1$\\\hline
 \end{tabular}
\end{center}

Note that $\widehat{\chi_{24}}$ is irreducible for every $\ell\neq
3$ by \cite[p. 140, 142, 143]{Atlas2}. We will show that
$\widehat{\chi_{24}}|_{3.P}$ is also irreducible for every $\ell\neq
3$. The structure of $P$ is $[3^5]:2S_4$. We denote by $O_3$ the
maximal normal $3$-subgroup (of order $3^5$) of $P$. Then the order
of any element in $O_3$ is $1$, $3$ or $9$, and $3.O_3$ is the
maximal normal $3$-subgroup of $3.P$. Since the values of
$\chi_{24}$ is zero at any element $\overline{g}$ where the order of
$g$ is $3$ or $9$, we have $$[\chi_{24}|_{3.O_3},
\chi_{24}|_{3.O_3}]_{3.O_3}=\frac{1}{|O_3|}\sum_{g\in
O_3}\chi_{24}(\overline{g})\overline{\chi_{24}(\overline{g})}=\frac{1}{3^5}27^2=3.$$
By Clifford's theorem,
$\chi_{24}|_{3.O_3}=e\cdot\sum_{i=1}^t\theta_i$ where
$e=[\chi_{24}|_{3.O_3}, \theta_1]_{3.O_3}$ and $\theta_1$,
$\theta_2$, ..., $\theta_t$ are the distinct conjugates of
$\theta_1$ in $3.P$. So $e^2t=3$ and therefore $e=1$, $t=3$. Thus,
$\chi_{24}|_{3.O_3}=\theta_1+\theta_2+\theta_3$. By Lemma
\ref{orbitreduction}, $\widehat{\chi_{24}}|_{3.P}$ is irreducible
when $\ell\neq 3$.

Arguing similarly, we also have $\widehat{\chi_{25}}|_{3.P}$ is
irreducible and therefore $\widehat{\overline{\chi_{24}}}|_{3.P}$ as
well as $\widehat{\overline{\chi_{25}}}|_{3.P}$ are irreducible for
every $\ell\neq 3$. Note that an outer automorphism $\tau$ of
$G_2(3)$ sends $P$ to $Q$ and fixes $\{\chi_{24}$,
$\overline{\chi_{24}}$, $\chi_{25}$, $\overline{\chi_{25}}\}$. Hence
the Lemma also holds when $M=3.Q$ by Lemma \ref{lemmatau}.
\end{proof}

\begin{lemma}\label{smalllemma4}
Theorem \ref{main theorem2} holds in the case $G=3\cdot G_2(3)$,
$M=3.(U_3(3):2)$.
\end{lemma}

\begin{proof}
Since the Schur multiplier of $U_3(3)$ is trivial and
$\mathbb{Z}_3=Z(3\cdot G_2(3))\leq Z(M)$, $3.(U_3(3):2)=3\times
(U_3(3):2)$. So
$\mathfrak{m}_\mathbb{C}(M)=\mathfrak{m}_\mathbb{C}(U_3(3):2)=64$.
Therefore if $\varphi|_M$ is irreducible then $\varphi(1)=27$ and
$\varphi$ is the reduction modulo $\ell\neq 3$ of any of the four
irreducible complex characters of degree $27$ which we denote by
$\chi_{24}$, $\chi_{25}$, $\overline{\chi_{24}}$ and
$\overline{\chi_{25}}$ as before.

When $\ell=2$ or $7$, $U_3(3):2$ does not have any irreducible
$\ell$-Brauer character of degree $27$. So we assume that $\ell\neq
2,3$ and $7$. That means $\ell\nmid |M|$ and
$\IBR_{\ell}(M)=\Irr(M)$. Therefore we only need to consider the
complex case.

It is obvious that $\varphi|_{3\times(U_3(3):2)}$ is irreducible if
and only if $\varphi|_{(U_3(3):2))}$ is irreducible. Moreover, since
$\varphi(1)=27$, $\varphi|_{U_3(3):2)}$ is irreducible if and only
if $\varphi|_{U_3(3)}$ is irreducible. Now we will show that either
$\chi_{24}|_{U_3(3)}$ or $\chi_{25}|_{U_3(3)}$ is irreducible and
the other is reducible.

From \cite[p. 14]{Atlas1}, we know that $U_3(3)$ has the unique
irreducible character of degree $27$, which is denoted by
$\chi_{10}$. It is easy to see that the classes $4A$, $4B$ of
$U_3(3)$ is contained in the same class of $G_2(3)$, which is $4A$
or $4B$. With no loss we suppose that this class is $4A$. Then the
class $4C$ of $U_3(3)$ is contained in the class $4B$ of $G_2(3)$.
In $U_3(3)$ we have $(8A)^2\subset 4A$, $(8B)^2\subset 4B$ and in
$G_2(3)$, $(8A)^2\subset 4A$, $(8B)^2\subset 4B$. So the classes
$8A$, $8B$ of $U_3(3)$ are contained in the class $8A$ of $G_2(3)$.
Similarly, the classes $12A$, $12B$ of $U_3(3)$ are contained in the
class $12A$ of $G_2(3)$. Comparing the values of
$\chi_{24}|_{U_3(3)}$ as well as $\chi_{25}|_{U_3(3)}$ with those of
$\chi_{10}$ on every conjugacy classes of $U_3(3)$, we see that
$\chi_{25}|_{U_3(3)}=\chi_{10}$ and $\chi_{24}|_{U_3(3)}$ is
reducible.

Note that $G_2(3)$ has two non-conjugate maximal subgroups which are
isomorphic to $U_3(3):2$. We denote these groups by $M_1$ and $M_2$.
Suppose that $\chi_{25}|_{M_1}$ is irreducible and
$\chi_{24}|_{M_1}$ is reducible. Let $\tau$ be an automorphism of
$G_2(3)$ such that $\tau(M_1)=M_2$. Then $\tau(M_2)=M_1$,
$\chi_{24}\circ\tau=\chi_{25}$ and $\chi_{25}\circ\tau=\chi_{24}$.
By Lemma \ref{lemmatau}, $\chi_{25}|_{M_2}$ is reducible and
$\chi_{24}|_{M_2}$ is irreducible.
\end{proof}

\begin{lemma}\label{smalllemma5}
Theorem \ref{main theorem2} holds in the case $G=3\cdot G_2(3)$,
$M=3.(L_3(3):2)$ or $3.(L_2(8):3)$.
\end{lemma}

\begin{proof}
This Lemma is proved similarly as the previous Lemma.
\end{proof}

\begin{lemma}\label{smalllemma7}
Theorem \ref{main theorem2} holds in the case $G=G_2(4)$,
$M=U_3(4):2$.
\end{lemma}

\begin{proof}
We have $\mathfrak{m}_{\mathbb{C}}(U_3(4))=75$ and
$\mathfrak{m}_{\mathbb{C}}(U_3(4):2)=150$. So if $\varphi|_M$ is
irreducible then $\varphi(1)\leq 150$. Inspecting the character
table of $G_2(4)$, we see that $\varphi(1)=64$, $78$ when $\ell=3$
or $\varphi(1)=65$, $78$ when $\ell\neq 2$, $3$. Using Lemma
\ref{lemma} and arguing as in Case $4$ of the proof of Theorem
\ref{main theorem}, we see that the restriction of the character of
smallest degree of $G_2(4)$ to $U_3(4)$ is irreducible.

It remains to consider the case that $\varphi$ is the reduction
modulo $\ell\neq2$ of the character $\chi_3$ (as denoted in \cite[p.
98]{Atlas1}) of degree $78$. Because $U_3(4):2$ has no complex
irreducible character of degree $78$, $\chi_3|_{U_3(4):2}$ is
reducible and so is $\widehat{\chi_3}|_{U_3(4):2}$.
\end{proof}

\begin{lemma}\label{smalllemma8}
Theorem \ref{main theorem2} holds in the case $G=G_2(4)$, $M=J_2$.
\end{lemma}

\begin{proof}
Comparing the degrees of irreducible $\ell$-Brauer characters of
$G_2(4)$ with those of $J_2$, we see that if $\varphi|_{J_2}$ is
irreducible then $\varphi(1)$ only can be any of the two irreducible
characters of degree $300$ with $\ell\neq 2$, $3$ and $7$. When
$\ell=0$, these two complex characters are denoted by $\chi_4$ and
$\chi_5$ in \cite[p. $98$]{Atlas1}. First, we will show that
$\chi_4|_{J_2}$ is actually irreducible. More precisely,
$\chi_4|_{J_2}=\chi_{20}$, where $\chi_{20}$ is the unique
irreducible complex character of $J_2$ of degree $300$.

It is easy to see that the values of $\chi_4$ and $\chi_{20}$ are
the same at conjugacy classes of elements of order $5$, $7$, $8$,
$10$, and $15$. The unique class $12A$ of elements of order $12$ in
$J_2$ is real. So it is contained in a real class of $G_2(4)$.
Therefore it is contained in class $12A$ of $G_2(4)$ and we have
$\chi_4(12A)=\chi_{20}(12A)=1$. We see that $(12A)^3=4A$ in both
$G_2(4)$ and $J_2$. Therefore the class $4A$ of $J_2$ is contained
in the class $4A$ of $G_2(4)$ and we also have
$\chi_4(4A)=\chi_{20}(4A)=4$. Now we consider classes of elements of
order $2$, $3$ and $6$. Since $J_2$ is a subgroup of $G_2(4)$,
either $2\cdot J_2$ (the universal cover of $J_2$) or $2\times J_2$
is a subgroup of $2\cdot G_2(4)$. Note that
$\mathfrak{d}_{\mathbb{C}}(2\cdot G_2(4))=12$ and
$\mathfrak{d}_{\mathbb{C}}(J_2)=\mathfrak{d}_{\mathbb{C}}(2\times
J_2)=14$. So $2\times J_2$ cannot be a subgroup of $2\cdot G_2(4)$
and therefore $2\cdot J_2$ is a subgroup of $G_2(4)$. From
\cite{Atlas1}, the class $2A$ of $G_2(4)$ lifts to two involution
classes of $2\cdot G_2(4)$ and the class $2B$ of $G_2(4)$ lifts to a
class of elements of order $4$ of $2\cdot G_2(4)$. In the same way,
the class $2A$ of $J_2$ lifts to two involution classes of $2\cdot
J_2$ and the class $2B$ of $J_2$ lifts to a class of elements of
order $4$ of $2\cdot J_2$. These imply that the classes $2A$ and
$2B$ of $J_2$ are contained in the classes $2A$ and $2B$ of
$G_2(4)$, respectively. Again, we have
$\chi_4(2A)=\chi_{20}(2A)=-20$ and $\chi_4(2B)=\chi_{20}(2B)=0$.
Using similar arguments, we also can show that the classes $6A$ and
$6B$ of $J_2$ are contained in the classes $6A$ and $6B$ of
$G_2(4)$, respectively. Again, $\chi_4(6A)=\chi_{20}(6A)=1$ and
$\chi_4(6B)=\chi_{20}(6B)=0$. In both $G_2(4)$ and $J_2$, we have
$(6A)^2=3A$ and $(6B)^2=3B$. That means the classes $3A$ and $3B$ of
$J_2$ are contained in the classes of $3A$ and $3B$ of $G_2(4)$,
respectively. One more time, $\chi_4(3A)=\chi_{20}(3A)=-15$ and
$\chi_4(3B)=\chi_{20}(3B)=0$. We have shown that the values of
$\chi_4$ and $\chi_{20}$ agree at all conjugacy classes of $J_2$.
Therefore $\chi_4|_{J_2}=\chi_{20}$. Note that
$\chi_5=\overline{\chi_4}$ and all irreducible characters of $J_2$
are real. Therefore we also have $\chi_5|_{J_2}=\chi_{20}$.

When $\ell\neq 2$, $3$, and $7$, the reductions modulo $\ell$ of
$\chi_4$, $\chi_5$ and $\chi_{20}$ are still irreducible. Thus,
$\widehat{\chi_4}|_{J_2}$ and $\widehat{\chi_5}|_{J_2}$ are
irreducible for every $\ell\neq 2$, $3$, and $7$.
\end{proof}

\begin{lemma}\label{smalllemma9}
Theorem \ref{main theorem2} holds in the case $G=2\cdot G_2(4)$,
$M=2.P$.
\end{lemma}

\begin{proof}
We have $\mathfrak{m}_\mathbb{C}(2.P)\leq
\sqrt{|2.P/Z(2.P)|}\leq\sqrt{|P|}=\sqrt{184,320}<430$. Therefore if
$\varphi|_M$ is irreducible then $\varphi(1) < 430$. There are five
cases to consider:

$\bullet$ $\varphi(1)=12$ when $\ell\neq 2$. Then $\varphi$ is the
reduction modulo $\ell\neq 2$ of the unique irreducible complex
character of degree $12$ of $2\cdot G_2(4)$. Throughout the proof of
this lemma we denote this character by $\chi$. Now we will show that
$\chi|_{2.P}$ is irreducible.

Note that if $g_1$ and $g_2$ are the pre-images of an element $g\in
G_2(4)$ under the natural projection $\pi: 2\cdot G_2(4)\rightarrow
G_2(4)$, then $\chi(g_1)=\pm\chi(g_2)$. Therefore $[\chi|_{2.P},
\chi|_{2.P}]_{2.P}=\frac{1}{2\cdot |P|}\sum_{x\in 2.P}
\chi(x)\overline{\chi(x)}=\frac{1}{|P|}\sum_{g\in
P}\chi(\overline{g})\overline{\chi(\overline{g})}$, where
$\overline{g}$ is a pre-image of $g$. We choose $\overline{g}$ so
that the value of $\chi$ at $\overline{g}$ is printed in \cite[p.
98]{Atlas1}. In \cite[p. 357]{EY}, we have the fusion of conjugacy
classes of $P$ in $G_2(4)$. By comparing the orders of centralizers
of conjugacy classes of $G_2(4)$ in \cite[p. 364]{EY} with those in
\cite[p. 98]{Atlas1} and looking at the values of irreducible
characters of degrees $65$, $78$, we can find a correspondence
between conjugacy classes of $G_2(4)$ in these two papers. The
length of each conjugacy class of $P$ can be computed from \cite[p.
357]{EY}. All the above information is collected in Table V. From
this table, we get
\begin{center}
\begin{tabular}{lll}$\sum_{g\in
P}\chi(\overline{g})\overline{\chi(\overline{g})}$&$=$&$12^2+3\cdot(-4)^2+
60\cdot(-4)^2+120\cdot(-4)^2+1440\cdot(-4)^2$\\
&&$+5760\cdot 2^2+320\cdot(-6)^2+960\cdot2^2+3840\cdot2^2+3840\cdot2^2$\\
&&$+3840\cdot 2^2+2\cdot3072\cdot(-3)^2+2\cdot9216\cdot1^2$\\
&$=$&$184,320$.
\end{tabular}
\end{center}
Therefore, $[\chi|_{2.P}, \chi|_{2.P}]_{2.P}=\frac{1}{|P|}\sum_{g\in
P}\chi(\overline{g})\overline{\chi(\overline{g})}=1$, which implies
that $\chi|_{2.P}$ is irreducible.

Next, we show that $\widehat{\chi}|_{2.P}$ is also irreducible for
every $\ell\neq 2$. Set $O_2=2^{2+8}$ to be the maximal normal
$2$-subgroup of $P$. Then $O_2$ is a union of conjugacy classes of
$P$. Note that $O_2$ is a $2$-group of exponent $4$ and so the
orders of elements in $O_2$ are $1$, $2$ or $4$. Hence $O_2$ is
either $A_0\cup A_1\cup A_2 \cup A_3 \cup A_{41}\cup A_{42}\cup A_5$
or $A_0\cup A_1\cup A_2 \cup A_{61}$. If the latter case happens
then $[\chi|_{2.O_2}, \chi|_{2.O_2}]=\frac{1}{|O_2|}\sum_{g\in
O_2}\chi(\overline{g})\overline{\chi(\overline{g})}=
(12^2+3\cdot(-4)^2+60\cdot(-4)^2)/1024=1152/1024$ which is not an
integer. Therefore the former case must happen and we have
$[\chi|_{2.O_2}, \chi|_{2.O_2}]_{2.O_2}=\frac{1}{|O_2|}\sum_{g\in
O_2}\chi(\overline{g})\overline{\chi(\overline{g})}=
(12^2+3\cdot(-4)^2+60\cdot(-4)^2+120\cdot(-4)^2)/1024=3$. By
Clifford's theorem, $\chi|_{2.O_2}=e\cdot\sum_{i=1}^{t}\theta_i$,
where $e=[\chi|_{2.O_2}, \theta_1]_{2.O_2}$ and $\theta_1$,
$\theta_2$, ..., $\theta_t$ are the distinct conjugates of
$\theta_1$ in $2.P$. So $e^2t=3$ and therefore $e=1$, $t=3$. Thus
$\chi|_{2.O_2}=\theta_1+\theta_2+\theta_3$. By Lemma
\ref{orbitreduction}, $\widehat{\chi}|_{2.P}$ is irreducible for
every $\ell\neq 2$.

\medskip
\centerline{{\sc {Table V.}} Fusion of conjugacy classes of $P$ in
$G_2(4)$}

\begin{center}
\begin{tabular}{|c|c|c|c|} \hline
Fusion in \cite{EY}& Corresponding class \cite{Atlas1}&Length&
 Value of $\chi$\\\hline
 $A_0\subset A_0$&$1A$& $1$&$12$\\\hline
 $A_1\subset A_1$&$2A$& $3$&$-4$\\\hline
 $A_2\subset A_1$&$2A$&$60$&$-4$\\\hline
 $A_3\subset A_2$&$2B$&$240$&$0$\\\hline
 $A_{41}\subset A_{31}$&$4A$&$120$&$-4$\\\hline
 $A_{42}\subset A_{32}$&$4C$&$360$&$0$\\\hline
 $A_{5}\subset A_{4}$&$4B$&$240$&$0$\\\hline
 $A_{61}\subset A_{2}$&$2B$&$960$&$0$\\\hline
 $A_{62}\subset A_{31}$&$4A$&$1440$&$-4$\\\hline
 $A_{63}\subset A_{32}$&$4C$&$1440$&$0$\\\hline
 $A_{71}\subset A_{51}$&$8A$&$5760$&$0$\\\hline
 $A_{72}\subset A_{52}$&$8B$&$5760$&$2$\\\hline
 $B_{0}\subset B_0$&$3A$&$320$&$-6$\\ \hline
 $B_{1}\subset B_1$&$6A$&$960$&$2$\\ \hline
 $[B_{2}]\subset B_1$&$6A$&$3840$&$2$\\ \hline
 $[B_{3}]\subset B_1$&$6A$&$3840$&$2$\\ \hline
 $B_{2}(0)\subset B_2(0)$&$12A$&$3840$&$2$\\ \hline
 $B_{2}(i)\subset B_2(i)$ $(i=1, 2)$
 &$12B, 12C$&$3840$&$0$\\\hline
 $C_{31}(i)\subset C_{21}$ (two classes)
 &$3B$&$5120$&$0$\\\hline
$C_{32}(i)\subset C_{22}$ (two classes)&$6B$&$15360$&$0$\\\hline
 $C_{41}(i)\subset C_{21}$ (two classes)&$3B$&$1024$&$0$\\\hline
 $C_{42}(i)\subset C_{22}$ (two classes)&$6B$&$15360$&$0$\\\hline
 $D_{11}(i)\subset D_{11}(i)$ (two classes)&$5C, 5D$&$3072$&$-3$\\\hline
 $D_{12}(i)\subset D_{12}(i)$ (two classes)&$10A, 10B$&$9216$&$1$\\\hline
 $E(i)\subset E_{1}(i)$ (four classes)&$15C, 15D$&$12288$&$0$\\\hline
\end{tabular}
\end{center}

$\bullet$ $\varphi(1)=104$ when $\ell\neq 2$, $5$. In this case,
$\varphi$ is actually the reduction modulo $\ell$ of a complex
irreducible character of degree $104$. Note that $104\nmid
368,640=|2.P|$ and therefore this case is not an example.

$\bullet$ $\varphi(1)=364$ when $\ell\neq 2$, $3$. Then $\varphi$ is
the reduction modulo $\ell$ of a unique faithful irreducible complex
character of degree $364$. Again, $364\nmid 368,640=|2.P|$.

$\bullet$ $\varphi(1)=352$ when $\ell=3$. Denote by $\chi_{364}$ the
unique faithful irreducible complex character of degree $364$ of
$2\cdot G_2(4)$. Then $\varphi=\widehat{\chi_{364}}-\widehat{\chi}$.
Suppose that $\varphi|_{2.P}$ is irreducible. Because
$\chi_{364}|_{2.P}$ is reducible and $\widehat{\chi}|_{2.P}$ is
irreducible, $\chi_{364}|_{2.P}=\lambda+\mu$ where $\lambda$ and
$\mu$ are irreducible complex characters of $2.P$ such that
$\widehat{\lambda}=\widehat{\chi}|_{2.P}$ and $\widehat{\mu}\in
\IBR_{3}(2.P)$. So $\mu\in \Irr(2.P)$ and $\mu(1)=352$ which leads
to a contradiction since $352\nmid |2.P|$.

$\bullet$ $\varphi(1)=92$ when $\ell=5$. Denote by $\chi_{104}$ one
of the two complex irreducible characters of degree $104$ of $2\cdot
G_2(4)$. Then $\varphi=\widehat{\chi_{104}}-\widehat{\chi}$. Suppose
that $\varphi|_{2.P}$ is irreducible. Because $\chi_{104}|_{2.P}$ is
reducible and $\widehat{\chi}|_{2.P}$ is irreducible,
$\chi_{104}|_{2.P}=\lambda+\mu$ where $\lambda$ and $\mu$ are
irreducible complex characters of $2.P$ such that
$\widehat{\lambda}=\widehat{\chi}|_{2.P}$ and $\widehat{\mu}\in
\IBR_{5}(2.P)$. So $\mu\in \Irr(2.P)$ and $\mu(1)=92$. Again, we get
a contradiction since $92\nmid |2.P|$.
\end{proof}

\begin{lemma}\label{smalllemma10}
Theorem \ref{main theorem2} holds in the case $G=2\cdot G_2(4)$,
$M=2.Q$.
\end{lemma}

\begin{proof}
By similar arguments as in Lemma \ref{smalllemma9}, we only need to
show that $\widehat{\chi}|_{2.Q}$ is irreducible for every
$\ell\neq2$, where $\chi$ is the unique irreducible character of
degree $12$ of $G$. First, let us prove that $\chi|_{2.Q}$ is
irreducible.

We have $[\chi|_{2.Q}, \chi|_{2.Q}]_{2.Q}=\frac{1}{2\cdot
|Q|}\sum_{x\in 2.Q}
\chi(x)\overline{\chi(x)}=\frac{1}{|Q|}\sum_{g\in
Q}\chi(\overline{g})\overline{\chi(\overline{g})}$, where
$\overline{g}$ is a pre-image of $g$ under $\pi$. We choose
$\overline{g}$ so that the value of $\chi$ at $\overline{g}$ is
printed in \cite[p. 98]{Atlas1}. In \cite[p. 361]{EY}, we have the
fusion of conjugacy classes of $Q$ in $G_2(4)$ and the length of
each conjugacy class of $Q$, which are described in Table VI. From
this table, we get
\begin{center}
\begin{tabular}{lll}$\sum_{g\in
Q}\chi(\overline{g})\overline{\chi(\overline{g})}$&$=$&$12^2+15\cdot(-4)^2+
360\cdot(-4)^2+240\cdot(-4)^2+240\cdot(-4)^2$\\
&&$+5760\cdot 2^2+2\cdot64\cdot(-6)^2+2\cdot960\cdot2^2+2\cdot3840\cdot2^2$\\
&&$+2\cdot3840\cdot2^2+2\cdot3072\cdot 2^2+4\cdot12288\cdot(-1)^2$\\
&$=$&$184,320$.
\end{tabular}
\end{center}
Therefore, $[\chi|_{2.Q}, \chi|_{2.Q}]_{2.Q}=\frac{1}{|Q|}\sum_{g\in
Q}\chi(\overline{g})\overline{\chi(\overline{g})}=1$, which implies
that $\chi|_{2.Q}$ is irreducible.

Next, we show that $\widehat{\chi}|_{2.Q}$ is also irreducible. Set
$O_2=2^{4+6}$ to be the maximal normal $2$-subgroup of $Q$. Since
$\chi|_{2.Q}$ is irreducible, by Clifford's theorem,
$\chi|_{2.O_2}=e\cdot\sum_{i=1}^{t}\theta_i$, where
$e=[\chi|_{2.O_2}, \theta_1]_{2.O_2}$ and $\theta_1$, $\theta_2$,
..., $\theta_t$ are the distinct conjugates of $\theta_1$ in $2.Q$.
We have $12=\chi(1)=et\theta_1(1)$. Note that $O_2$ is a $2$-group
of exponent $4$ and so $\theta_1(1)\in\{1, 2, 4\}$. Therefore $et\in
\{3, 6, 12\}$. Set $m=e^2t=[\chi|_{2.O_2},
\chi|_{2.O_2}]_{2.O_2}=\frac{1}{|O_2|}\sum_{g\in
O_2}\chi(\overline{g})\overline{\chi(\overline{g})}$. From the above
table, we see that $\chi(\overline{g})=0$ or $-4$ for every element
$g$ of order $2$ or $4$. Therefore $m=\frac{1}{|O_2|}\sum_{g\in
O_2}\chi(\overline{g})\overline{\chi(\overline{g})}=\frac{1}{1024}
(12^2+n\cdot(-4)^2)$ where $n$ is the sum of the lengths of
conjugacy classes of $O_2$ at which the values of $\chi$ is $-4$.
Then $n=(1024\cdot m-144)/16$. We also have
$n\leq15+360+240+240=855$ and $5\mid n$. This implies that $m\leq
13$. Because $et\in\{3, 6, 12\}$, now we see that $m=e^2t\in \{3, 6,
9, 12\}$. If $m=3$, resp. $9$, $12$, then $n=183$, resp. $556$,
$759$, which are coprime to $5$. So $m=6$ and then $n=375$ which is
the sum of lengths of classes $A_1$ and $A_{32}$. We have shown that
$e^2t=6$ and therefore $e=1$, $t=6$. Therefore
$\chi|_{2.O_2}=\sum_{i=1}^{6}\theta_i$, which implies that
$\widehat{\chi}|_{2.Q}$ is irreducible for every $\ell\neq 2$ by
Lemma \ref{orbitreduction}.
\end{proof}

\medskip
\centerline{{\sc {Table VI.}} Fusion of conjugacy classes of $Q$ in
$G_2(4)$}

\begin{center}
\begin{tabular}{|c|c|c|c|} \hline
Fusion in \cite{EY}& Corresponding class \cite{Atlas1}&Length&
 Value of $\chi$\\\hline
 $A_0\subset A_0$&$1A$& $1$&$12$\\\hline
 $A_1\subset A_1$&$2A$& $15$&$-4$\\\hline
 $A_2\subset A_2$&$2B$&$48$&$0$\\\hline
 $A_{31}\subset A_2$&$2B$&$240$&$0$\\\hline
 $A_{32}\subset A_{31}$&$4A$&$360$&$-4$\\\hline
 $A_{33}\subset A_{32}$&$4C$&$360$&$0$\\\hline
 $A_{41}\subset A_{1}$&$2A$&$240$&$-4$\\\hline
 $A_{42}(0)\subset A_{31}$&$4A$&$240$&$-4$\\\hline
 $A_{42}(i)\subset A_{4}$ $(i=1, 2)$&$4B$&$240$&$0$\\\hline
 \begin{tabular}{c}$A_{5}(t)\subset A_{2}, A_{32}, A_4$\\
 (four classes)
 \end{tabular}&$2B, 4C, 4B$&$720$&$0$\\\hline
 $A_{61}\subset A_{51}$&$8A$&$5760$&$0$\\\hline
 $A_{62}\subset A_{52}$&$8B$&$5760$&$2$\\\hline
 $B_{0}(i)\subset B_0$ $(i=1, 2)$&$3A$&$64$&$-6$\\ \hline
  $B_{1}(i)\subset B_1$ $(i=1, 2)$&$6A$&$960$&$2$\\ \hline
$B_{2}(i)\subset B_1$ $(i=1, 2)$&$6A$&$3840$&$2$\\\hline
  $B_{3}(i,0)\subset B_2(0)$ $(i=1, 2)$&$12A$&$3840$&$2$\\ \hline
\begin{tabular}{c}  $B_{3}(i,j)\subset B_2(j)$\\
$(i=1, 2; j=1, 2)$\end{tabular}&$12B, 12C$&$3840$&$0$\\ \hline
 $C_{21}\subset C_{21}$&$3B$&$5120$&$0$\\\hline
$C_{22}\subset C_{22}$&$6B$&$15360$&$0$\\\hline
 $C_{31}(i)\subset C_{21}$ (two classes)&$3B$&$5120$&$0$\\\hline
 $C_{32}(i)\subset C_{22}$ (two classes)&$6B$&$15360$&$0$\\\hline
 $D_{11}(i)\subset D_{21}(i)$ (two classes)&$5A, 5B$&$3072$&$2$\\\hline
 $D_{12}(i)\subset D_{22}(i)$ (two classes)&$10C, 10D$&$9216$&$0$\\\hline
 $E(i)\subset E_{2}(i)$ (four classes)&$15A, 15B$&$12288$&$-1$\\\hline
\end{tabular}
\end{center}

\begin{lemma}\label{smalllemma11}
Theorem \ref{main theorem2} holds in the case $G=2\cdot G_2(4)$,
$M=2.(U_3(4):2)$ or $M=2.(SL_3(4):2)$.
\end{lemma}

\begin{proof}
We give the details of the proof in the case $M=2.(U_3(4):2)$. Note
that the Schur multiplier of $U_3(4)$ is trivial. So
$M=2.(U_3(4):2)=(2\times U_3(4)).2$. Therefore
$\mathfrak{m}_{\mathbb{C}}(M)\leq 2\mathfrak{m}_{\mathbb{C}}(2\times
U_3(4))=150$ by \cite[p. 30]{Atlas1}. Hence, if $\varphi|_M$ is
irreducible then $\varphi(1)\leq 150$. From the character table of
$2\cdot G_2(4)$, we have $\varphi(1)=12$ when $\ell\neq 2$,
$\varphi(1)=104$ when $\ell\neq 2, 5$ or $\varphi(1)=92$ when
$\ell=5$.

$\bullet$ $\varphi(1)=92$ when $\ell=5$. Inspecting the $5$-Brauer
character table of $U_3(4)$ in \cite[p. 72]{Atlas2}, we see that $M$
does not have any irreducible $5$-Brauer character of degree $92$.
So this case is not an example.

$\bullet$ $\varphi(1)=12$ when $\ell\neq 2$. In this case, $\varphi$
is the reduction modulo $\ell\neq2$ of the unique irreducible
complex character $\chi$ of degree $12$ of $2\cdot G_2(4)$.
Inspecting the character tables of $U_3(4)$, we have
$\mathfrak{d}_{\ell}(U_3(4))=12$ for every $\ell\neq 2$. It follows
that $\widehat{\chi}|_{U_3(4)}$ is irreducible and so is
$\widehat{\chi}|_M$.

$\bullet$ $\varphi(1)=104$ when $\ell\neq 2$, $5$. In this case,
$\varphi$ is the reduction modulo $\ell$ of any of the two faithful
irreducible complex characters of degree $104$ of $2\cdot G_2(4)$,
which are denoted by $\chi_{34}$ and $\chi_{35}$ as in \cite[p.
98]{Atlas1}. First, we will show that $\chi_{34}|_M$ is irreducible.

From Lemma \ref{smalllemma7}, we know that the restriction of the
unique irreducible character $\chi_2$ of degree $65$ of $G_2(4)$ to
$U_3(4)$ is irreducible and equal to the unique rational irreducible
character of degree $65$ of $U_3(4)$. By looking at the values of
these characters, we see that the involution class $2A$ of $U_3(4)$
is contained in the class $2A$ of $G_2(4)$ and the class $5E$ of
$U_3(4)$ is contained in the class $5A$ or $5B$ of $G_2(4)$. With no
loss we assume that this class is $5A$.

Now suppose that $\chi_{34}|_{U_3(4)}$ contains the unique
irreducible character of degree $64$ of $U_3(4)$. Since
$\chi_{34}(2A)=8$, by looking at the values of irreducible
characters of $U_3(4)$ at the involution class $2A$, we see that the
other irreducible constituents of $\chi_{34}|_{U_3(4)}$ are of
degrees $1$ and $39$. But then we get a contradiction by looking at
the values of these characters at the class $5E$ of $U_3(4)$.

So $\chi_{34}|_{U_3(4)}$ does not contain the irreducible character
of degree $64$ of $U_3(4)$. Inspecting the character table of
$U_3(4)$ in \cite[p. 30]{Atlas1}, we see that $U_3(4)$ has four
irreducible characters of degree $52$, which are denoted by
$\chi_9$, $\chi_{10}$, $\chi_{11}$ and $\chi_{12}$. Note that
$\chi_{10}=\overline{\chi_9}$ and $\chi_{12}=\overline{\chi_{11}}$.
We also see that the values of irreducible characters (of degrees
different from $64$) of $U_3(4)$ at the class $5E$ are: $1, 2, -b5,
b5+1, -b5+1, b5+2, b5, -b5-1, 0$ where
$b5=\frac{1}{2}(-1+\sqrt{5})$. Moreover, this value is $b5$ if and
only if the character is $\chi_9$ or $\chi_{10}$. Note that
$\chi_{34}(5A)=2b5$. Now we suppose that $2b_5=x_1\cdot 1+x_2\cdot
2+ x_3\cdot (-b_5)+x_4\cdot (b_5+1)+x_5\cdot (-b_5+1)+x_6\cdot
(b_5+2)+x_7\cdot b_5+x_8\cdot (-b_5-1)$ where $x_i$s $(i=1,2,...,8)$
are nonnegative integer numbers. Then $x_1+2x_2+x_4+x_5+2x_6-x_8=0$
and $-x_3+x_4-x_5+x_6+x_7-x_8=2$. These equations imply that
$2=-x_3+x_4-x_5+x_6+x_7-(x_1+2x_2+x_4+x_5+2x_6)=-x_1-2x_2-2x_5+x_7$.
Therefore $x_7\geq 2$. In other words, there are at least two
irreducible constituents of degree $52$ in $\chi_{34}|_{U_3(4)}$.
Note that $\chi_{34}(1)=104$. So $\chi_{34}|_{U_3(4)}$ must be the
sum of two irreducible characters of degree $52$. Since $\chi_{34}$
is real, it is easy to see that
$\chi_{34}|_{U_3(4)}=\chi_9+\chi_{10}$. Notice that $M=(2\times
U_3(4)).2$ and the ``2'' is an outer automorphism of $U_3(4)$ that
fuses $\chi_9$ and $\chi_{10}$. So $\chi_{34}|M$ is irreducible.

It is easy to see that $\chi_{35}=\ast\circ\chi_{34}$ and
$\chi_{11}+\chi_{12}=\ast\circ(\chi_9+\chi_{10})$ where the operator
$\ast$ is the algebraic conjugation: $r+s\sqrt{5}\mapsto
r-s\sqrt{5}$ for $r, s\in \mathbb{Q}$. Since
$\chi_{34}|_{U_3(4)}=\chi_{9}+\chi_{10}$,
$\chi_{35}|_{U_3(4)}=\chi_{11}+\chi_{12}$, which implies that
$\chi_{35}|_M$ is irreducible.

When $\ell\neq 2,5$, the reductions modulo $\ell$ of $\chi_{34}$ and
$\chi_{35}$ are still irreducible. Arguing similarly as above, we
have $\widehat{\chi_{34}}|_M$ as well as $\widehat{\chi_{35}}|_M$
are irreducible.
\end{proof}

{\bf Proof of Theorem \ref{main theorem2}}.

\medskip
(i) According to \cite[p. $61$]{Atlas1}, if $M$ is a maximal
subgroup of $G_2(3)$ then $M$ is $G_2(3)$-conjugate to one of the
following groups:
\begin{enumerate}
  \item[1)] $P=[3^5]:2S_4$, $Q=[3^5]:2S_4$, maximal parabolic subgroups,
  \item[2)] $U_3(3):2$, two non-conjugate subgroups,
  \item[3)] $L_3(3):2$, two non-conjugate subgroups,
  \item[4)] $L_2(8):3$,
  \item[5)] $2^3\cdot L_3(2)$,
  \item[6)] $L_2(13)$,
  \item[7)] $[2^5]:3^2.2$.
\end{enumerate}

By Lemmata $\ref{smalllemma1} - \ref{smalllemma2}$, we need to
consider the following cases:

1) $M=P$ or $Q$. Since the structure of $P$ as well as $Q$ is
$[3^5]:2S_4$, they are solvable. It is well-known that every Brauer
character of $M$ is liftable to a complex character. Therefore the
degree of every irreducible Brauer character of $M$ divides
$|P|=|Q|=11,664$ and less than $\sqrt{11,664}=108$. Checking both
the complex and Brauer character tables of $G_2(3)$, we see that
there is no character satisfying these conditions.

3) $M=L_3(3):2$. From \cite[p. 13]{Atlas1}, we have
$\mathfrak{m}_{\mathbb{C}}(L_3(3):2)=52$. So if $\varphi|_{M}$ is
irreducible then $\varphi(1)\leq 52$ and therefore $\varphi(1)=14$.
That means $\varphi$ is the reduction modulo $\ell\neq3$ of the
unique irreducible complex character of degree $14$, which is
denoted by $\chi_2$ in \cite[p. 60]{Atlas1}. Since $14\nmid
|L_3(3):2|=11,232$, $\chi_2|_{L_3(3):2}$ as well as
$\widehat{\chi_2}|_{L_3(3):2}$ are reducible for every $\ell\neq3$.

The cases $M=L_2(8):3, L_2(13)$, and $[2^5]:3^2.2$ are treated
similarly.

\medskip
(ii) In this part, we only consider faithful irreducible characters
of $3\cdot G_2(3)$. They are characters which are not inflated from
irreducible characters of $G_2(3)$. By Lemma \ref{maximalsubgroup},
a maximal subgroup of $3\cdot G_2(3)$ is the pre-image of a maximal
subgroup of $G_2(3)$ under the natural projection $\pi: 3\cdot
G_2(3)\rightarrow G_2(3)$. We denote by $3.X$ the pre-image of $X$
under $\pi$. By Lemmata $\ref{smalllemma3} - \ref{smalllemma5}$, we
need to consider the following cases:

5) $M=3.(2^3\cdot L_3(2))$. If $\varphi|_M$ is irreducible then
$\varphi(1)\leq\mathfrak{m}_\mathbb{C}(M)\leq\sqrt{|(2^3\cdot
L_3(2))|}=\sqrt{1344}<37$. So $\varphi(1)=27$ and $\varphi$ is the
reduction modulo $\ell\neq 3$ of an irreducible character $\chi$ of
degree $27$ of $3\cdot G_2(3)$. Since $27\nmid |M|$, $\chi|_M$ is
reducible and so is $\widehat{\chi}|_M$.

6) $M=3.L_2(13)$. From \cite[p. 8]{Atlas1}, the Schur multiplier of
$L_2(13)$ has order $2$. So $3.L_2(13)=3\times L_2(13)$. Hence
$\mathfrak{m}_\mathbb{C}(M)=\mathfrak{m}_\mathbb{C}(L_2(13))=14$. On
the other hand, the degree of any faithful irreducible Brauer
character of $3\cdot G_2(3)$ is at least $27$. So we do not have any
example in this case.

The case $M=3.([2^5]:3^2.2)$ is similar.

\medskip
(iii) According to \cite[p. $97$]{Atlas1}, if $M$ is a maximal
subgroup of $G_2(4)$ then $M$ is $G_2(4)$-conjugate to one of the
following groups:
\begin{enumerate}
  \item[1)] $P=2^{2+8}:(3\times A_5)$, $Q=2^{4+6}:(A_5\times 3)$, maximal parabolic
  subgroups,
  \item[2)] $U_3(4):2$,
  \item[3)] $SL_3(4):2$,
  \item[4)] $U_3(3):2$,
  \item[5)] $A_5\times A_5$,
  \item[6)] $L_2(13)$,
  \item[7)] $J_2$.
\end{enumerate}

By Lemmata $\ref{smalllemma7} - \ref{smalllemma8}$, we need to
consider the following cases:

1) $M=P$ or $Q$. From \cite{EY}, it is easy to see that
$\mathfrak{m}_{\mathbb{C}}(P)$ as well as
$\mathfrak{m}_{\mathbb{C}}(Q)$ are less than $256$. So if
$\varphi|_M$ is irreducible then $\varphi(1)<256$. Inspecting the
complex and Brauer character tables of $G_2(4)$, we have two
possibilities:

$\bullet$ $\ell\neq 2$, $3$ and $\varphi$ is the reduction modulo
$\ell$ of any of the two characters of degree $65$ and $78$, which
are $\chi_2$ and $\chi_3$ as denoted in \cite[p. 98]{Atlas1}. Note
that both $65$ and $78$ do not divide $|P|=|Q|=184,320$. So both
$\widehat{\chi_2}|_{P,Q}$ and $\widehat{\chi_3}|_{P,Q}$ are
reducible.

$\bullet$ $\ell=3$ and $\varphi(1)=64$ or $78$. The case
$\varphi(1)=78$ cannot happen by a similar reason as above. If
$\varphi(1)=64$ then
$\varphi=\widehat{\chi_2}-\widehat{1}_{G_2(4)}$. Suppose that the
restriction of this character to $P$ is irreducible. Since
$\chi_2|_P$ is reducible, $\chi_2|_{P}=\lambda+\mu$ where $\lambda$,
$\mu\in\Irr(P)$, $\widehat{\lambda}=\widehat{\mathbf{1}_{P}}$ and
$\widehat{\mu}\in\IBR_3(P)$. We then have $\mu(1)=64$. Inspecting
the character table of $P$ in \cite[p. 358]{EY}, we see that there
is no irreducible character of $P$ of degree $64$ and we get a
contradiction. The argument for $Q$ is exactly the same.

3) $M=SL_3(4):2$. This case is treated similarly as Case $3$ when
$q\equiv1(\bmod$ $3)$ in the proof of Theorem \ref{main theorem}.

The cases $M=U_3(3):2$, $A_5\times A_5$, and $L_2(13)$ are similar.

\medskip
(iv) By Lemma \ref{maximalsubgroup}, a maximal subgroup of $2\cdot
G_2(4)$ is the pre-image of a maximal subgroup of $G_2(4)$ under the
natural projection $\pi: 2\cdot G_2(4)\rightarrow G_2(4)$. We denote
by $2.X$ the pre-image of $X$ under $\pi$. By Lemmata
$\ref{smalllemma9} - \ref{smalllemma11}$, we need to consider the
following cases.

4) $M=2.(U_3(3):2)$. Since the Schur multiplier of $U_3(3)$ is
trivial, $2.U_3(3)=2\times U_3(3)$ and $2.(U_3(3):2)=(2\times
(U_3(3)).2$. So if $\varphi|_M$ is irreducible then $\varphi(1)\leq
2\mathfrak{m}_\mathbb{C}(U_3(3))=64$. Therefore $\varphi$ is the
reduction modulo $\ell\neq 2$ of $\chi$, the unique irreducible
complex character of degree $12$ of $2\cdot G_2(4)$. Assume that
$\chi|_M$ is irreducible. Using the character table of $U_3(3)$ in
\cite[p. 14]{Atlas1}, it is easy to see that
$\chi|_{U_3(3)}=2\chi_2$, where $\chi_2$ is the unique irreducible
character of degree $6$ of $U_3(3)$. Therefore $\chi|_{2\times
U_3(3)}=2(\sigma\otimes \chi_2)$, where $\sigma$ is the nontrivial
irreducible character of $\mathbb{Z}_2=Z(2\cdot G_2(4))$. This and
Lemma \ref{lemmaisaacs} imply that $\chi|_M$ is reducible, a
contradiction.

5) $M=2.(A_5\times A_5)$. We denote by $A$ and $B$ the pre-images of
the first and second terms $A_5$ (in $A_5\times A_5$), respectively,
under the projection $\pi$. For every $a\in A$, $b\in B$, we have
$\pi([a,b])=[\pi(a), \pi(b)]=1$. Therefore, $[a,b]\in \mathbb{Z}_2$
where $\mathbb{Z}_2=Z(2\cdot G_2(4))\leq Z(M)$. This implies that
$[[A, B], A]=[[B, A], A]=1$. By the $3$-subgroup lemma, we have
$[[A, A], B]=1$. Since the Schur multiplier of $A_5$ is $2$, $A$ is
$2\times A_5$ or $2.A_5$, the universal cover of $A_5$. If $A=2.A_5$
then $[A, A]=A$. If $A=2\times A_5$ then $[A, A]=A_5$. So, in any
case, $[A, B]=1$ or $A$ centralizes $B$. That means $M=(A\times
B)/\mathbb{Z}_2$.

We have
$\mathfrak{m}_\mathbb{C}(M)\leq\sqrt{|M/Z(M)|}\leq\sqrt{|A_5\times
A_5|}=60$. Therefore if $\varphi|_M$ is irreducible then $\varphi$
is the reduction modulo $\ell\neq 2$ of $\chi$, the unique complex
irreducible character of degree $12$ of $2\cdot G_2(4)$. Now we will
show that $\chi|_M$ is reducible.

Suppose $\lambda$ is any irreducible character of $M$ of degree
$12$. Then we can regard $\lambda$ as an irreducible character of
$A\times B$ with $\mathbb{Z}_2\subset \Ker \lambda$. Assume that
$\lambda=\lambda_A\otimes\lambda_B$ where $\lambda_A\in \Irr(A)$ and
$\lambda_B\in\Irr(B)$. There are two possibilities:

$\bullet$ One of $\lambda_A(1)$ and $\lambda_B(1)$ is $2$ and the
other is $6$. With no loss of generality, assume that
$\lambda_A(1)=2$ and $\lambda_B(1)=6$. From the character tables of
$A_5$ and $2.A_5$ in \cite[p. 2]{Atlas1}, we see that $A_5$ does not
have any irreducible character of degree $2$ or $6$. So $A=2.A_5$.
The value of any irreducible character of degree $2$ of $2.A_5$ at
any conjugacy class of elements of order $5$ is
$\frac{1}{2}(-1\pm\sqrt{5})$. Therefore, $\lambda|_A=6\lambda_A$ is
not rational.

$\bullet$ One of $\lambda_A(1)$ and $\lambda_B(1)$ is $4$ and the
other is $3$. With no loss of generality, assume that
$\lambda_A(1)=4$ and $\lambda_B(1)=3$. From the character tables of
$A_5$ and $2.A_5$ in \cite[p. 2]{Atlas1}, we see that the value of
any irreducible character of degree $3$ of $A_5$ or $2.A_5$ at any
conjugacy class of elements of order $5$ is
$\frac{1}{2}(1\pm\sqrt{5})$. Therefore, $\lambda|_B=4\lambda_B$ is
not rational.

We have shown that any irreducible character of degree $12$ of $M$
is not rational. On the other hand, $\chi$ is rational. Hence
$\chi|_M$ is reducible and so is $\widehat{\chi}|_M$.

The cases $M=2.L_2(13)$ and $2.J_2$ are treated similarly.
\hfill$\Box$

\section{Results for Suzuki and Ree groups}

\begin{theorem}
Let $G=Sz(q)$ be the Suzuki group where $q=2^n$, $n$ is odd and
$n\geq 3$. Let $\varphi$ be an absolutely irreducible character of
$G$ in characteristic $\ell\neq2$ and $M$ a maximal subgroup of $G$.
Assume that $\varphi(1)>1$. Then $\varphi|_M$ is irreducible if and
only if $M$ is $G$-conjugate to the maximal parabolic subgroup of
$G$ and $\varphi$ is the reduction modulo $\ell$ of any of the two
irreducible complex characters of degree $(q-1)\sqrt{q/2}$.
\end{theorem}

\begin{proof}
According to \cite{Su}, if $M$ is a maximal subgroup of $G$, then
$M$ is $G$-conjugate to one of the following groups:
\begin{enumerate}
  \item[1)] $P=[q^2].\mathbb{Z}_{q-1}$, the maximal parabolic subgroup,
  \item[2)] $D_{2(q-1)}$,
  \item[3)] $\mathbb{Z}_{q+\sqrt{2q}+1}.\mathbb{Z}_4$,
  \item[4)] $\mathbb{Z}_{q-\sqrt{2q}+1}.\mathbb{Z}_4$,
  \item[5)] $Sz(q_0)$, $q=q_0^\alpha$, $\alpha$ prime, $q_0\geq8$.
\end{enumerate}

By Lemmata \ref{1}, \ref{2} and the irreducibility of $\varphi|_M$,
we have $\sqrt{|M|}\geq \mathfrak{d}_{\ell}(G)$, which is larger or
equal to $(q-1)\sqrt{q/2}$ by \cite{T1}. Therefore, $|M|\geq
q(q-1)^2/2$. This inequality happens if and only if $M$ is the
maximal parabolic subgroup of $G$.

The complex character table of $P$ is given in \cite[p. 157]{M}.
From this table, we have
$\mathfrak{m}_{\mathbb{C}}(P)=(q-1)\sqrt{q/2}$. Since
$\mathfrak{d}_{\ell}(G)\leq \varphi(1)\leq\mathfrak{m}_{\ell}(P)\leq
\mathfrak{m}_{\mathbb{C}}(P)$, it follows that
$\varphi(1)=(q-1)\sqrt{q/2}$. Using the notation and results about
Brauer trees of the Suzuki group in \cite{B}, we have
$\varphi=\widehat{\Gamma_1}$ or $\widehat{\Gamma_2}$, where
$\Gamma_1$ and $\Gamma_2$ are the two irreducible complex characters
of $Sz(q)$ of degree $(q-1)\sqrt{q/2}$. Now we will show that the
restrictions of these characters to $P$ is indeed irreducible.

If $\ell =0$, comparing directly the values of characters on
conjugacy classes, we see that $\Gamma_1|_P=\phi_2$ and
$\Gamma_2|_P=\phi_3$, where $\phi_2$ and $\phi_3$ are irreducible
characters of $P$ and their values are given in \cite[p. 157]{M}.

Next, we will show that $\widehat{\phi_i}$s, $i=2$, $3$, are
irreducible when $\ell\neq 0,2$. Consider the element $f$ of order
$4$ in $P$ which is given in \cite[p. 157]{M}. Since $\ell$ is odd
and $\ord (f)=4$, $f$ is an $\ell$-regular element. Assume the
contrary that $\widehat{\phi_2}$ is reducible. Then it is the sum of
more than one irreducible Brauer characters of $P$ whose degrees are
less than $(q-1)\sqrt{q/2}$. Since $P$ is solvable, every Brauer
character of $P$ is liftable to complex characters. Inspecting the
complex character table of $P$, we see that the value at the element
$f$ of any irreducible character of degree less than
$(q-1)\sqrt{q/2}$ is real. On the other hand, $\phi_2(f)$ is not
real, a contradiction. We have shown that $\widehat{\phi_2}$ is
irreducible and so is $\widehat{\phi_3}$, as
$\phi_3=\overline{\phi_2}$.
\end{proof}

{\bf{Note}}: Most of the Schur multipliers of the Suzuki group and
the Ree group are trivial except the Schur multiplier of $Sz(8)$,
which is an elementary abelian group of order $2^2$. Note that
$2^2=Z(2^2.Sz(8))$ which is not cyclic. Therefore, $2^2.Sz(8)$ does
not have any faithful irreducible character. This multiplier $2^2$
has three cyclic quotients of order $2$ which are corresponding to
groups $2.Sz(8)$, $2'.Sz(8)$ and $2''.Sz(8)$. These groups are
permuted by the automorphism group. Let us consider our problem for
one of them, say $2.Sz(8)$.

Inspecting the character tables of $2.Sz(8)$, we see that if
$\varphi$ is a faithful irreducible character of $2.Sz(8)$, then
$\varphi(1)\geq8$. Therefore, if $\varphi|_M$ is irreducible then
$\sqrt{|M/Z(M)|}\geq8$. So the unique possibility for $M$ is $2.P$,
where $P$ is the maximal parabolic subgroup of $Sz(8)$. Moreover,
$\varphi(1)\leq\mathfrak{m}_\mathbb{C}(2.P)\leq\sqrt{|P|}=\sqrt{448}<22$.
Hence, we have $\varphi(1)=8$ when $\ell=5$ or $\varphi(1)=16$ when
$\ell=13$.

$\bullet$ If $\varphi(1)=8$ when $\ell=5$ then
$\varphi=\varphi_{11}$ as denoted in \cite[p. 64]{Atlas2}. Note that
$P=[2^6].7$ and $2.P=2.([2^6].7)=[2^7].7$. From \cite[p.
64]{Atlas2}, the value of $\varphi_{11}$ at any nontrivial
$2$-element is $0$. Therefore, $[\varphi_{11}|_{[2^7]},
\varphi_{11}|_{[2^7]}]_{[2^7]}=8^2/2^6=1$ and hence
$\varphi_{11}|_{[2^7]}$ is irreducible. So $\varphi_{11}|_{2.P}$ is
also irreducible.

$\bullet$ If $\varphi(1)=16$ when $\ell=13$ then
$\varphi=\varphi_{9}$ as denoted in \cite[p. 65]{Atlas2}. Since
$16>\sqrt{2^7}$, $\varphi_{9}|_{[2^7]}$ is reducible. By Lemma
\ref{lemmaisaacs}, $\varphi_{9}|_{2.P}$ is also reducible.

So when $G=2.Sz(8)$, if $\varphi$ is faithful, $\varphi|_M$ is
irreducible if and only if $M=2.P$ and $\varphi$ is the unique
irreducible $5$-Brauer character of degree $8$.

\begin{theorem}
Let $G=\ta G_2(q)$ be the Ree group where $q=3^n$, $n$ is odd and
$n\geq 3$. Let $\varphi$ be an absolutely irreducible character of
$G$ in characteristic $\ell\neq3$ and $M$ be a maximal subgroup of
$G$. Assume that $\varphi(1)>1$. Then $\varphi|_M$ is irreducible if
and only if $M$ is $G$-conjugate to the maximal parabolic subgroup
of $G$ and $\varphi$ is the nontrivial constituent (of degree
$q^2-q$) of the reduction modulo $\ell=2$ of the unique irreducible
complex character of degree $q^2-q+1$.
\end{theorem}

\begin{proof}
According to \cite[p. 181]{K2}, if $M$ is a maximal subgroup of $G$,
then $M$ is $G$-conjugate to one of the following groups:
\begin{enumerate}
  \item[1)] $P=[q^3]:\mathbb{Z}_{q-1}$, the maximal parabolic subgroup,
  \item[2)] $2\times L_2(q)$, involution centralizer,
  \item[3)] $(2^2\times D_{(q+1)/2}):3$,
  \item[4)] $\mathbb{Z}_{q+\sqrt{3q}+1}:\mathbb{Z}_6$,
  \item[5)] $\mathbb{Z}_{q-\sqrt{3q}+1}:\mathbb{Z}_6$,
  \item[6)] $^2G_2(q_0)$, $q=q_0^\alpha$, $\alpha$ prime.
\end{enumerate}

By Lemmata \ref{1}, \ref{2} and the irreducibility of $\varphi|_M$,
we have $\sqrt{|M|}\geq \mathfrak{d}_{\ell}(G)$, which is larger or
equal to $q(q-1)$ by \cite{T1}. Therefore, $|M|\geq q^2(q-1)^2$.
This inequality happens if and only if $M$ is the maximal parabolic
subgroup $P$. The complex character table of $P$ is given in
\cite[p. 88]{LM}. From there, we get $\mathfrak{m}_{\mathbb{C}}(P)=
q(q-1)$ and therefore $\mathfrak{m}_{\ell}(P)\leq q(q-1)$.

Assume that $\ell=0$ or $\ell\geq 5$. Using the results about Brauer
trees of $G$ in \cite{H3}, it is easy to check that
$\mathfrak{d}_{\ell}(G)=q^2-q+1>\mathfrak{m}_{\ell}(P)$, which
contradicts Lemma \ref{1}. So it remains to consider $\ell=2$.

We have $q(q-1)\geq\mathfrak{m}_{\ell}(P)\geq \varphi(1)\geq
\mathfrak{d}_{\ell}(G)\geq q(q-1)$. Therefore $\varphi(1)=q(q-1)$.
We will check all 2-blocks of $G$ which are studied in \cite{W} and
\cite{LM}. We also use the notation in these papers.

1) $\xi_9$, $\xi_{10}$, $\eta_i^{\pm}$ are of 2-defect $0$. Their
degrees are all larger than $q(q-1)$.

2) There is one 2-block of defect 1. All characters in this block
are $\eta_r$ and $\eta'_r$ whose degrees are $q^3+1$. By Lemma
\ref{Navaro}, there is a unique irreducible 2-Brauer character of
degree $q^3+1$ in this block.

3) There are several 2-blocks of defect $2$. Every character in
these blocks has degree $(q-1)(q^2-q+1)$. Applying Lemma
\ref{Navaro} again, all irreducible 2-Brauer characters in these
blocks have degree $(q-1)(q^2-q+1)$.

4) The principal block and its decomposition matrix is described in
\cite{LM}. We have $\varphi_1(1)=1$, $\varphi_2(1)=q(q-1)$,
$\varphi_3(1)=(q-1)(q^2-\sqrt{\frac{q}{3}}(q+1)+1)> q(q-1)$,
$\varphi_4(1)=\varphi_5(1)=(q-1)\sqrt{\frac{q}{3}}(q+1-3\sqrt{\frac{q}{3}})/2
>q(q-1)$.

In summary, the unique possibility for $\varphi$ is
$\varphi=\varphi_2=\widehat{\xi_2}-\widehat{1_G}$ when $\ell=2$,
where $\xi_2$ is the unique irreducible complex character of degree
$q^2-q+1$ of $G$. Now we will prove that $\varphi_2|_P$ is indeed
irreducible. Suppose that $\varphi_2|_P$ is reducible. Then it is
the sum of more than one irreducible $2$-Brauer characters of $P$.
These characters are of degrees less than $q(q-1)$. Moreover, since
$P$ is solvable, they are liftable. Therefore, their values at the
element $X$ of order $3$ (which is a representative of a conjugacy
class of $P$ given in \cite[p. 88]{LM} is positive. On the other
hand, $\varphi_2(X)=-q$ which is negative, a contradiction.
\end{proof}


\begin{thebibliography}{ABCD}
\bibitem[A]{A}
  M. Aschbacher, On the maximal subgroups of the finite classical
  groups, {\it Invent. Math.} {\bf 76} $(1984)$, $469-514$.

\bibitem[AH]{AH}
  J. An and S. Huang, The character tables of the parabolic subgroups of the Chevalley
  groups of type $G_2(q)$, {\it Comm. Algebra} {\bf 34} $(2006)$, $1763-1792$.

\bibitem[Atlas1]{Atlas1}
  J. H. Conway, R. T. Curtis, S. P. Norton, R. A. Parker, and R. A. Wilson,
 {\it `ATLAS of Finite Groups'}, Clarendon Press, Oxford,
$1985$.

\bibitem[Atlas2]{Atlas2}
  C. Jansen, K. Lux, R. Parker, R. Wilson, {\it `An ATLAS of Brauer
  characters'}, Clarendon Press, Oxford, $1995$.

\bibitem[B]{B}
  R. Burkhardt, {\" U}ber die Zerlegungszahlen der Suzukigruppen
  $Sz(q)$, {\it J. Algebra} {\bf 59} $(1979)$, $421-433$.

\bibitem[Cp]{Cp}
  B. N. Cooperstein, Maximal subgroups of $G_2(2^n)$, {\it J. Algebra}
   {\bf 70} $(1981)$, $23-36$.

\bibitem[CR]{CR}
  B. Chang and R. Ree, The characters of $G_2(q)$,
  {\it Symposia Mathematica, Vol. 13, Academic Press, London}
   $(1974)$, $395-413$.

\bibitem[E]{E}
  H. Enomoto, The characters of $G_2(q)$, $q= 3^n$, {\it Japan. J. Math.}
   {\bf 2} $(1976)$, $191-248$.

\bibitem[EY]{EY}
  H. Enomoto and H. Yamada, The characters of $G_2(2^n)$, {\it Japan. J. Math.}
   {\bf 12} $(1986)$, $325-377$.

\bibitem[F]{F}
 W. Feit, {\it `The representation theory of finite groups'},
 North-Holland Publ. Comp., Amsterdam, New York, Oxford, $1982$.

\bibitem[G]{G}
 M. Geck, Irreducible Brauer characters of the 3-dimensional
 special unitary groups in non-defining characteristic, {\it Comm. Algebra}
   {\bf 18} $(1990)$, $563-584$.

\bibitem[GT]{GT}
  R. M. Guralnick and P. H. Tiep,
  Low-dimensional representations of special linear groups in cross characteristic,
   {\it Proc. London Math. Soc.} {\bf 78} $(1999)$, $116-138$.

\bibitem[H1]{H1}
  G. Hiss, On the decomposition numbers of $G_2(q)$, {\it J.
  Algebra} {\bf 120} $(1989)$, $339-360$.

\bibitem[H2]{H2}
  G. Hiss, Zerlegungszahlen endlicher Gruppen vom Lie-Typ
  in nicht-definierender Charakteristik, Habilitationsschrift, RWTH Aachen, 1990.

\bibitem[H3]{H3}
  G. Hiss, The Brauer trees of the Ree groups, {\it Comm. Algebra}
  {\bf 19} $(1991)$, $871-888$.

\bibitem[HS1]{HS1}
  G. Hiss and J. Shamash, 3-blocks and 3-modular characters of $G_2(q)$,
{\it J. Algebra} {\bf 131} $(1990)$, $371 - 387$.

\bibitem[HS2]{HS2}
  G. Hiss and J. Shamash, 2-blocks and 2-modular characters of $G_2(q)$,
{\it Math. Comp.} {\bf 59} $(1992)$, $645 - 672$.

\bibitem[I] {I}
  I. M. Isaacs, {\it `Character theory of finite groups'}, Dover
  Publications, Inc., New York, $1994$.

\bibitem[K1]{K1}
 P. B. Kleidman, The maximal subgroups of the Chevalley groups $G_2(q)$ with $q$ odd,
 the Ree groups $^2G_2(q)$, and their automorphism groups, {\it J. Algebra}
   {\bf 117} $(1988)$, $30-71$.

\bibitem[K2]{K2}
  P. B. Kleidman, The subgroup structure of some finite simple
  groups, Ph. D thesis, Trinity College, $1986$.

\bibitem[KT]{KT}
  A. I. Kostrikin and P. H. Tiep, {\it `Orthogonal Decompositions and Integral
  Lattices'}, Walter de Gruyter, Berlin a.o., $1994$.

\bibitem[LM]{LM}
  P. Landrock and G. O. Michler, Principle $2$-blocks of the
  simple groups of Ree type, {\it Trans. Amer. Math. Soc.} {\bf
  260} $(1980)$, $83-111$.

\bibitem[M]{M}
  R. P. Martineau, On representation of the Suzuki groups over
  fields of odd characteristic, {\it J. London Math. Soc.} {\bf 6}
  $(1972)$, $153-160$.

\bibitem[OT]{OT}
  T. Okuyama and Y. Tsushima, Local properties of $p$-block algebras
  of finite groups, {\it Osaka J. Math.} {\bf 20} $(1983)$, $33-41$.

\bibitem[S]{S}
  J. Saxl, unpublished manuscript.

\bibitem[S1]{S1}
  J. Shamash, Blocks and Brauer trees for groups of type $G_2(q)$,
  {\it The Arcata Conference on Representations of Finite Groups,
   Part 2}, Amer. Math. Soc., Providence, RI, $(1987)$, $283 - 295$.

\bibitem[S2]{S2}
  J. Shamash, Brauer trees for blocks of cyclic defect in the
  groups $G_2(q)$ for primes dividing $q^2\pm q+1$, {\it J.
  Algebra} {\bf 123} $(1989)$, $378-396$.

\bibitem[S3]{S3}
  J. Shamash, Blocks and Brauer trees in the
  groups $G_2(q)$ for primes dividing $q\pm 1$, {\it Comm. Algebra
  } {\bf 17} $(1989)$, $1901-1949$.

\bibitem[S4]{S4}
  J. Shamash, Blocks and Brauer trees for the
  groups $G_2(2^k)$, $G_2(3^k)$, {\it Comm.
  Algebra} {\bf 20} $(1992)$, $1375-1387$.

\bibitem[SF]{SF}
  W. A. Simpson and J. S. Frame, The character tables for $SL(3,q)$,
  $SU(3,q)$, $PSL(3,q)$, $PSU(3,q)$, {\it Can. J. Math.}, Vol. XXV, No.
  3 $(1973)$, $486-494$.

\bibitem[Su]{Su}
  M. Suzuki, On a class of doubly transitive groups, {\it Ann. of
  Math.} {\bf 75} $(1962)$, $105-145$.

\bibitem[T1]{T1}
  P. H. Tiep, Low dimensional representations of finite quasisimple groups,
{\it Groups, combinatorics and geometry} (Durham, 2001), World Sci.
Publishing, River Edge, NJ, (2003), 277--294.

\bibitem[T2]{T2}
  P. H. Tiep, Finite groups admitting grassmannian $4$-designs, {\it J. Algebra}
  {\bf 306} $(2006)$, $227-243$.

\bibitem[TZ]{TZ}
  P. H. Tiep and A. E. Zalesskii, Minimal polynomials of semisimple elements of
  finite classical groups in cross characteristic representations, in preparation.

\bibitem[W]{W}
 H. N. Ward, On Ree's series of simple groups, {\it Trans. Amer. Math. Soc.}
   {\bf 121} $(1966)$, $62-89$.
\end{thebibliography}
\end{document}